\documentclass[onecolumn]{autart}
\usepackage{}    

\usepackage{graphicx}           
\usepackage{dcolumn}            
\usepackage{bm}                 

\usepackage{graphics}           
\usepackage{epsfig}             
\usepackage{times}              
\usepackage[fleqn]{amsmath}
\usepackage{amssymb}
\usepackage{extarrows}
\usepackage{amsfonts}
\usepackage{mathrsfs}
\usepackage{fancyhdr}
\usepackage{color}
\usepackage{subfigure}
\usepackage{enumerate}

\allowdisplaybreaks[4] 

\newtheorem{theorem}{Theorem}
\newtheorem{lemma}[theorem]{Lemma}
\newtheorem{corollary}[theorem]{Corollary}

\newtheorem{remark}{Remark}

\begin{document}

\begin{frontmatter}

\title{A minimum problem with free boundary and subcritical growth in Orlicz spaces}
\author{Jun Zheng $^{1}$}\ead{zhengjun@swjtu.edu.cn},
\author{Leandro S. Tavares$^{2}$}\ead{lean.mat.ufca@gmail.com},
\author{Claudianor O. Alves$^{3}$}\ead{coalves@mat.ufcg.br}
\address{$^{1}${School of Mathematics, Southwest Jiaotong University,
        Chengdu 611756, Sichuan, China}\\
        $^{2}$Centro de Ci\^{e}ncias e Tecnologia, Universidade
Federal do Cariri, Avenida Tenente Raimundo Rocha, CEP 58429-900, Juazeiro
do Norte, CE, Brazil\\
$^{3}${Unidade Acad\^{e}mica de Matem\'{a}tica, Universidade Federal de Campina Grande, CEP 58429-900, Campina Grande, PB, Brazil}}
\begin{abstract}
The aim of this paper is to study the heterogeneous optimization problem
 \begin{align*}
\mathcal {J}(u)=\int_{\Omega}(G(|\nabla
u|)+qF(u^+)+hu+\lambda_{+}\chi_{\{u>0\}} )\text{d}x\rightarrow\text{min},
\end{align*}
in the class of functions $ W^{1,G}(\Omega)$ with $ u-\varphi\in W^{1,G}_{0}(\Omega)$,
  for a given function $\varphi$, where
  $W^{1,G}(\Omega)$ is the class of weakly differentiable functions
  with $\int_{\Omega}G(|\nabla
u|)\text{d}x<\infty$. The functions $G$ and $F$ satisfy structural conditions of Lieberman's type that allow for a
 different behavior at $0$ and at $\infty$. {}{Moreover, $F$ allows for a subcritical growth.} Given functions $q,h$ and constant $\lambda_+\geq 0$, we address several regularity results for minimizers of $\mathcal {J}(u)$, including local $C^{1,\alpha}-$, and local Log-Lipschitz continuities for minimizers of $\mathcal {J}(u)$ with $\lambda_+=0$, and {}{$\lambda_+\geq 0$} respectively. We also establish growth rate near the free boundary for each non-negative minimizer of $\mathcal {J}(u)$ with $\lambda_+=0$, and $\lambda_+>0$ respectively. Furthermore, under additional assumption that $F\in C^1([0,+\infty); [0,+\infty))$, local Lipschitz regularity is carried out for non-negative minimizers of $\mathcal {J}(u)$ with $\lambda_{+}>0$.
\end{abstract}
\end{frontmatter}
\section{Introduction}\label{Sec: introduction}
  Let $\Omega$ be a smooth bounded domain in $\mathbb{R}^{N}(N> 2
)$. Let $g,F\in C([0,\infty);[0,\infty))\cap C^1((0,\infty);(0,\infty))$ with $g(0)=F(0)=0$ satisfying the Lieberman's conditions, which were introduced in \cite{L} for a large class of degenerate/singular elliptic equations, i.e.,
\begin{align}\label{+031}
0<{\delta_0}\leq \frac{tg'(t)}{g(t)}\leq g_{0},\ \forall\  t>0,
\end{align}
with $1+\delta_0 < N$  and
\begin{align}\label{hip-F}
     1+\theta_0\leq \frac{F'(t)t}{F(t)}\leq 1+f_0,\ \forall\  t>0.
\end{align}
with $\theta_0,f_0$ satisfying {}{$0<1+\theta_0\leq 1+f_0 \leq \displaystyle\frac{N(1+\delta_0)}{N-(1+\delta_0)}$}.

The aim of this paper is to derive interior regularity
estimates for {}{the minimizers} of a large class of heterogeneous non-differentiable functionals
\begin{align}\label{main equation}
\mathcal {J}(u)=\int_{\Omega}(G(|\nabla
u|)+qF(u^+)+hu+\lambda_{+}\chi_{\{u>0\}} )\text{d}x\rightarrow\text{min},
\end{align}
among competing functions  $u\in \{u\in
L^{1}(\Omega):\int_{\Omega}G(|\nabla u|)\text{d}x<\infty,u=\varphi \ \text{on}\ \partial \Omega\}$, where $G(t)=\displaystyle\int_{0}^tg(s)\text{d}s$, {}{$\varphi\in L^{\infty}(\Omega)\cap W^{1,G}(\Omega)$}, $q\in L^{\infty}(\Omega)$ with $q \not\equiv 0$, $u^{+}= \max\{u,0\},{}{ h \in L^{m}(\Omega)(m\geq N)}$ and ${\lambda_+\geq 0}$ is a constant. Note that if no restriction is made on the sign of $h$, problem \eqref{main equation} has a minimizer that may change its sign near the free boundary $\partial\{u>0\}$. Therefore problem \eqref{main equation} is not in the one-phase case in the strict sense.

 A typical form of \eqref{main equation} is the free boundary problem of $p-$Laplacian, i.e.,
\begin{align}\label{main equation'}
\int_{\Omega}(|\nabla
u|^p+q(u^+)^\gamma+hu+\lambda_{+}\chi_{\{u>0\}} )\text{d}x\rightarrow\text{min},
\end{align}
over the set $\{u\in
W^{1,p}(\Omega),u-\varphi\
\in
W^{1,p}_0(\Omega)\}$, corresponding to set $g(t)=pt^{p-1},\delta_0=g_{0}=p-1,p>1$, and $F(t)= t^{\gamma},\theta_0=f_{0}=\gamma-1,0<\gamma<p$ in \eqref{+031}, and \eqref{hip-F} respectively {}{with $h$ being a measurable function.} More examples of functions satisfying \eqref{+031} (or \eqref{hip-F}) are given in Remark \ref{Remark 1}.

A number of important mathematical physics problems, coming from several different contexts,
are modeled by optimization setups, for which \eqref{main equation'} serves as an emblematic, leading
prototype. The case of $\gamma=1,q\not\equiv0$ and $ \lambda_{+}=0$ represents the obstacle type problems. The general case of $0<\gamma<p,q\not\equiv0$ and $ \lambda_{+}=0$ is usually used to model the density of certain chemical specie in reaction with a porous catalyst pellet. The case of $q\equiv 0$ and $ \lambda_{+}\neq 0$ relates to jets flow and cavities problems. The minimization problem \eqref{main equation'}, particularly homogeneous one-phase problem (i.e., $h\equiv0$, and minimizers of which do not change sign), has indeed received overwhelming attention at aspects of both regularity of solutions and regularity of free boundaries in the past decades, e.g., just to cite a few, \cite{KKPS,LS} for the homogeneous one-phase obstacle problems, \cite{AP,GG,P1,P2} for the homogeneous one-phase chemical reaction problems with $\gamma \in (0,1)$, \cite{ACF,DP} for the homogeneous one-phase jets flow and cavities problems with $q\equiv 0$ and $ \lambda_{+}\neq 0$, and \cite{LT} for a large class of homogeneous one-phase free boundary problems of $p-$Laplacian type corresponding to \eqref{main equation'} with $h\equiv0,1\leq \gamma <p$. {}{  We also point out the interesting work \cite{LdT} where the authors studied a two-phase version of \eqref{main equation'} with $0<\gamma\leq 1$ given by
$$
\int_{\Omega}(|\nabla
u|^p+(u^+)^\gamma+ (u^-)^\gamma + hu)\text{d}x\rightarrow\text{min},
$$
over the set $\{u\in
W^{1,p}(\Omega),u-\varphi\
\in
W^{1,p}_0(\Omega)\}$
 and $h \in L^{q}(\Omega)$ with $q \geq n.$}

 For the setting in Orlicz spaces, regularities of solutions and free boundaries are addressed for $\lambda_{+}=0,h\equiv 0, q\equiv C$ in \cite{CL,CLR}, and for $q\equiv h\equiv 0,\lambda_{+}>0$ in \cite{MW}, respectively. The homogeneous two-phase jets flow and cavities problems were studied in \cite{BM}. It should be mentioned that the heterogeneous two-phase problems related to \eqref{main equation} with $F(t)=t^{\gamma} (\gamma \in (0,1])$, and a version of two-phase problems related to \eqref{main equation} with $F(t)\leq \max\{t^{p},t^{p} \} (p,q\geq 1)$ and $h\equiv 0$, were studied in \cite{ZFZ,ZZZ}, and \cite{B} respectively. Nevertheless, regularity in problem \eqref{main equation} for a large class
of heterogeneous non-differentiable functionals has been little studied in the literature in Orlicz spaces.

The aim of this paper is to consider the free boundary problem \eqref{main equation} and prove several {}{regularity results} for minimizers of $\mathcal {J}(u)$. Comparing with the existing results, the main contribution of this paper include: {}{(i) It is proved the existence and boundness of minimizers in the critical case for problem \eqref{main equation'}. Moreover, an uniform $L^{\infty}$  estimate is proved for the subcritical case. At least to our knowledge such results are new in the literature even in the Laplacian case.} (ii) {}{With subcritical growth of $F$, several regularity results for minimizers are proved in the setting of Orlicz spaces, which are also new in the literature.} (iii) We establish local Log-Lipschitz continuity for minimizers of $\mathcal {J}(u)$ with $\lambda+>0$ under the assumption that $\delta_0>0$, which is weaker than the condition that $\delta_0\geq 1$ (or, equivalently, $\frac{g(t)}{t}$ is increasing in $t$) in \cite{B,ZZZ}; (iv) Our problems concern with not only the non-homogeneous case of $h\not\equiv 0$, but also the case of $F(t)\leq \max\{t^{p},t^{p} \} $ with positive exponents $p,q$ less or larger than $1$, which can be seen as complements of \cite{B,BM,ZFZ,ZZZ}; (v) We prove the growth rate near the free boundaries for non-negative minimizers of $\mathcal {J}(u)$ with $h\not\equiv 0$ and $F$ satisfying \eqref{hip-F}, which is new even in the problem \eqref{main equation'} with one-phase and $\gamma \in (0,p)$; (vi) we prove local Lipschitz continuity for non-negative minimizers of $\mathcal {J}(u)$ with $h\not\equiv 0$, which is an extension of \cite{B}.
{}{An important point regarding the results of this manuscript is that they are the  first steps in the understanding of analytic and geometric properties of the free boundary related to the  minimization problem \eqref{main equation'}.}

Throughout this paper, without spacial states, we always assume that
\begin{align*}
&g,{}{F}\in C([0,\infty);[0,\infty))\cap C^1((0,\infty);(0,\infty)),g(0)=F(0)=0,\ g,{}{F}\ \text{satisfy}\ \eqref{+031} \ \text{and}\  \eqref{hip-F},\\
&q\in L^{\infty}(\Omega),q \not\equiv 0, {}{h \in L^{m}(\Omega)(m\geq N)}, \varphi \in W^{1,G}(\Omega)\cap L^{\infty}(\Omega),\\
&\lambda_{+}\geq 0 \ \text{is\ a \ constant},
\end{align*}
where the definition of $W^{1,G}(\Omega)$ is given in Section \ref{Sec: auxiliary results}.
{}{For $t>0$, denote by $f(t)$ the derivative of $F(t)$, i.e, $f(t)= F'(t)$,  $\forall t>0$.} Let $\mathcal {K}=\{v \in W^{1,G}(\Omega): v-\varphi\in W^{1,G}_0(\Omega)\}$. Denote a ball in $\Omega$ by $B,$ $B_r$ or $B_R$ without special statements on their radius and centres, and denote by $B_r(x_0)$ a ball with radius $r$ and centre $x_0$. Without confusion, constants $\varepsilon,\tau,c,C,C_0,C_1,...$ appearing in this paper may be different from each other.

The rest of the paper is organized as follows. We provide two remarks on the structural conditions \eqref{+031} and \eqref{hip-F} at the end of Section \ref{Sec: introduction}. Some basic properties of functions $g,G$ and $F$, definitions of Orlicz spaces, properties of functions in Orlicz spaces, and an iteration lemma for the establishment of regularities of minimizers are presented in Section \ref{Sec: auxiliary results}. {}{With the critical growth of $F$,} existence of minimizers (and non-negative minimizers) of $\mathcal {J}(u)$ and their $L^\infty-$boundedness and local $C^{0,\alpha}-$continuity are addressed in Section \ref{Sec: existence and continuity}. {}{With the subcritical growth of $F$,} local $C^{1,\alpha}-$continuity, and local Log-Lipschitz continuity of minimizers are established in Section \ref{Sec: 4} for $\mathcal {J}(u)$ with $\lambda+=0$, and {}{$\lambda+\geq 0$} respectively. {}{With the subcritical growth of $F$,} growth rate near the free boundary $\partial\{u>0\}$ of each non-negative minimizer of $\mathcal {J}(u)$ with $\lambda+=0$ and $\lambda+>0$ are given respectively in Section \ref{Sec: 5}. As a consequence of the obtained results, we can prove the optimal growth rates of each non-negative minimizer and its gradient in the one-phase free boundary problems for $p-$Laplacian. Under the further assumption on $F$, i.e., $F\in C^1([0,+\infty); [0,+\infty))$, Local Lipschitz continuity of non-negative minimizers of $\mathcal {J}(u)$ with $\lambda+>0$ is established in Section \ref{Sec: Lipschitz continuity}.

\begin{remark}We do not require any $C^2-$continuity of $F$ to assume that
\begin{align}\label{hip-F'}
     \theta_0\leq \frac{f'(t)t}{f(t)}\leq f_0,\ \ \ \ \forall\  t>0,
\end{align}
provided a $C^1-$continuos function $f$ satisfying $f(0)=0$ and $F(t)=\int_{0}^tf(s)\text{d}s$. Therefore \eqref{hip-F} is weaker than the structural condition imposed on $F$ by \eqref{hip-F'}.
\end{remark}

\begin{remark} \label{Remark 1}
In this remark, we present several functions defined on $[0,+\infty)$ and satisfying condition \eqref{+031} {}{(or \eqref{hip-F})}, proofs of which and more functions satisfying a slight version of \eqref{+031} {}{(or \eqref{hip-F})} can be found in \cite{ZG}. {}{From these examples, one may find that the class of nonlinearities $F,$ caracterized by \eqref{hip-F}, is larger than  the ones considered in \cite{B,BM,LdT,ZFZ,ZZZ}.}
\begin{enumerate}[(i)]
\item $g(t)=(1+t)\ln(1+t)-t$ satisfies \eqref{+031} with $\delta_0=1$ and $g_0=2$.
\item $g(t)=\ln(1+at)+bt$ satisfies \eqref{+031} with $\delta_0=\frac{b}{a+b}$ and $g_0=\frac{a+b}{b}$, where $a>0,b>0$.
\item $g(t)=t^a\log_c(bt+d)$ satisfies \eqref{+031} with $\delta_0=a $ and $g_0=a+\frac{1}{\ln d}$, where $,a,b>0,c,d>1$.
\item $g(t)=\frac{t^a}{\log_c(bt+d)}$ satisfies \eqref{+031} with $\delta_0=a-\frac{1}{\ln d}$ and $g_0=a$, where $b>0,c,d>1, a>\frac{1}{\ln d}$.
\item $g(t)=\left\{\begin{aligned}
& at^p,0\leq t<t_0\\
&bt^q+c,t\geq t_0
\end{aligned}\right.$ satisfies \eqref{+031} with $\delta_0=\min\{p,q\}>0$ and $g_0=\max\{p,q\}$, where $a,b,c,p,q,t_0>0$ such that $  at_0^p=bt_0^q+c,$ and $apt_0^{p-1}=bqt_0^{q-1}.$
\end{enumerate}
\end{remark}
\section{Some auxiliary results}\label{Sec: auxiliary results}
\begin{lemma}[\cite{LU}]\label{lady-aux}
Let $(J_n), n \in \mathbb{N}$ a sequence of nonnegative numbers satisfying
$$J_{n+1} \leq CD^{n}J^{1+\xi}_{n}, n=0,1,..., $$
where $C,D>0$ and $\xi >1$ are constants that does not depend on $n \in \mathbb{N}.$ If
$$ J_0 \leq C^{-1\frac{1}{\xi}} D^{-\frac{1}{\xi^2}},$$
then $J_n \rightarrow 0$ as $n \rightarrow +\infty.$
\end{lemma}

\begin{lemma}[\cite{BM,MW}]\label{G-properties}The functions $g$ and $G$ satisfy the following
properties:
\begin{enumerate}
\item[$(g_{1})$] $\min\{s^{{\delta_0}}, s^{g_{0}} \}g(t) \leq  g(st) \leq
\max\{s^{{\delta_0}}, s^{g_{0}}\}g(t),\ \forall s,t\geq0.$
\item[$(G_{1})$] $G \text{\ is\  convex on } [0,+\infty)\text{\ and }\ C^{2}-\text{continuous on}\  (0,+\infty).$
 \item[$(G_{2})$]$\frac{tg(t)}{1+g_{0}}\leq G(t)\leq tg(t), \ \forall \ t\geq
 0.$
\item[$(G_{3})$]$\min\{s^{{\delta_0}+1}, s^{g_{0}+1} \}\frac{G(t)}{1+g_{0}}
\leq G(st) \leq
(1+g_{0})\max\{s^{{\delta_0}+1}, s^{g_{0}+1}\}G(t),\ \forall t\geq 0.$
\item[$(G_{4})$] $G(a+b)\leq 2^{g_{0}}(1+g_{0})(G(a)+G(b)), \ \forall \
 a,b\geq
 0.$
 \item[$(G_{5})$] $\delta_0\leq \frac{tg_s'(t)}{g_s(t)}\leq g_0$ and $\frac{1}{1+g_0}\leq G_s(1)\leq 1$ for all $t>0$, where $G_s(t)=\frac{G(st)}{sg(s)}$ and $ g_s(t)=G_s'(t)$ for $s>0$.
\end{enumerate}
\end{lemma}
\begin{lemma}\label{g-properties}{The functions $F$ and $f$ satisfy the following
properties:}
\begin{enumerate}
\item[$(F_{1})$]
 ${\min\{s^{1+\theta_0}, s^{1+f_0} \} F(t)\leq  F(st) \leq
\max\{s^{1+\theta_0}, s^{1+f_0}\}F(t)},\ \forall s,t\geq 0.$
\item[$(F_{2})$]${F(s+t)\leq 2^{1+f_0}(F(s)+F(t))},\ \forall s,t\geq 0.$
\item[$(f_{1})$]${\frac{1+\theta_0} {1+f_0}\min\{s^{\theta_0}, s^{f_0}\}f(t)\leq f(st)\leq \frac{1+f_0}{1+\theta_0}\max\{s^{\theta_0}, s^{f_0}\}f(t)},\ \forall s,t>0.$
 \item[$(f_2)$]
    {$
\lim\limits_{t\rightarrow +\infty} \frac{f(t)}{g(t)}=0$, and $
\lim\limits_{t\rightarrow +\infty} \frac{F(t)}{G(t)}=0$.
}
\end{enumerate}
\end{lemma}
\begin{pf*}{Proof}{$(F_{1})$ is a consequence of \eqref{hip-F} and $(g_1)$. }

{For $(F_{2})$, without loss of generality, assume that $s\leq t$, then $F(s+t)\leq F(2t)\leq 2^{1+f_0}F(t)\leq 2^{1+f_0}(F(s)+F(t))$.}

{For $(f_{1})$, we deduce by \eqref{hip-F} and $(F_1)$,
\begin{align*}
stf(st)\leq & (1+f_0)F(st)\\
\leq & (1+f_0)\max\{s^{1+\theta_0}, s^{1+f_0}\}F(t)\\
\leq & \frac{1+f_0}{1+\theta_0} \max\{s^{1+\theta_0}, s^{1+f_0}\}tf(t),
\end{align*}
  which yields the second inequality in $(f_{1})$. The first inequality in $(f_{1})$ can be obtained in a similar way.}

  {For $(f_{2})$, we infer from $(f_1)$ and $(g_1)$ that for large $t>1$,
\begin{align*}
f(t)t^{\delta_0-f_0}\leq \frac{1+f_0} {1+\theta_0} f(1)t^{f_0} t^{\delta_0-f_0}\leq \frac{1+f_0} {1+\theta_0} \frac{f(1)} {g(1)}   g(1)t^{\delta_0}\leq \frac{1+\theta_0} {1+f_0} \frac{f(1)} {g(1)}   g(t),
\end{align*}
which implies $ \lim\limits_{t\rightarrow +\infty} \frac{f(t)}{g(t)}=0$. The second result can be obtained in a similar way by $(F_1)$ and $(G_1)$.}
   \hfill $\blacksquare$
\end{pf*}
\begin{lemma}\label{f-properties} The follows statements hold true.
{\begin{enumerate}
\item[$(F_{3})$] If $f(t)$ is decreasing in $t>0$, then $F(s)-F(t)\leq F(s-t),\ \forall s\geq t\geq 0.$
 \item[$(F_{4})$] If $f(t)$ is increasing in $t>0$, then $F(s)-F(t)\leq f(M)(s-t),\ \forall 0\leq t\leq s\leq M$ with some $M>0$.
\end{enumerate}}
\end{lemma}
\begin{pf*}{Proof}
{We prove $(F_{3})$. Firstly, note that $F(0)=0$ due to $(F_1)$. Now Fix $s\geq 0$ and let $v(t)=F(s+t)- F(s)-F(t)$ for any $t\geq 0$. For $t>0$, we have $v'(t)=f(s+t)-f(t)\leq 0$, which yields the nondecreasing monotonicity of $v$ in $t>0$. By continuity of $v$ in $t=0$, we conclude that $v(t)\leq v(0)$ for all $t\geq 0$, i.e.,  $F(s+t)\leq F(s)+F(t)$. Finally, for $s\geq t\geq 0$, we have $F(s)-F(t)=F(s-t+t)-F(t)\leq F(s-t)+F(t)-F(t)=F(s-t)$.}

{For $(F_{4})$, if $s=t=0$, the result is trivial. If $M\geq s\geq t>0$, by Mean Value Theorem, it follows $F(s)-F(t)=f(\xi)(s-t)$ with some $\xi\in (t,s)\subset (0,M]$. By the increasing monotonicity of $f(t)$ in $t>0$, we have $F(s)-F(t)\leq f(M)(s-t)$ for all $M\geq s\geq t>0$. Finally, $(f_{4})$ has been proved.}
   \hfill $\blacksquare$
\end{pf*}
 As $g$ is strictly increasing we can define its inverse function $g^{-1}$. Then $g^{-1}$ satisfies a similar condition to (1.2).
\begin{lemma}[\cite{MW}]\label{Lemma 3}\ \ The function $g^{-1}$ satisfies the
following property:
\begin{align}
\frac{1}{g_{0}}\leq \frac{t(g^{-1})'(t)}{g^{-1}(t)}\leq
\frac{1}{{\delta_0}},\ \ \forall\ t>0.\notag
\end{align}
Moreover, $g^{-1}$ satisfies
\begin{align}
(\tilde{g}_{1})\ \ \min\{s^{\frac{1}{{\delta_0}}},
s^{\frac{1}{g_{0}}} \}g^{-1}(t) \leq g^{-1}(st) \leq
\max\{s^{\frac{1}{{\delta_0}}}, s^{\frac{1}{g_{0}}}\}g^{-1}(t),\notag
\end{align}
and if $\tilde{G}$ is such that $\tilde{G}'(t)=g^{-1}(t)$
then
\begin{align}
&(\tilde{G}_{1})\ \ \frac{1+{\delta_0}}{{\delta_0}}\min\{s^{1+1/{\delta_0}},
s^{1+1/g_{0}} \}\tilde{G}(t)\leq \tilde{G}(st) \leq
\frac{{\delta_0}}{1+{\delta_0}}\max\{s^{1+1/{\delta_0}},
s^{1+1/g_{0}}\}\tilde{G}(t),\notag\\
&(\tilde{G}_{2})\ \  ab\leq \varepsilon G(a)+C(\varepsilon)\tilde{G}(b),\ \
\forall\ a,b>0\ \text{and}\ \varepsilon>0,\notag\\
&(\tilde{G}_{3})\ \ \tilde{G}(g(t))\leq g_{0}G(t).\notag
\end{align}
\end{lemma}
We recall that the functional
\begin{align*}
\| u\|_{L^{G}(\Omega)}:= \inf \left\{k > 0; \int_{\Omega} G\left(\frac{|u|}{k}\right) \text{d}x \leq 1 \right\}
\end{align*}
is a norm in the Orlicz space $L^{G}(\Omega)$ which is the linear hull of the Orlicz class
\begin{align*}
\mathcal{K}_{G}(\Omega):= \left\{u \ \text{ measurable};  \int_{\Omega}G(|u|)\text{d}x   < +\infty  \right\}.
\end{align*}
Notice that this set is convex, since $G$ is also convex. The Orlicz-Sobolev space $W^{1,G}(\Omega)$ is defined as
\begin{align*}
W^{1,G}(\Omega):= \{u \in L^{G}(\Omega) ; \nabla u \in (L^{G}(\Omega))^n\},
\end{align*}
which is the usual subspace of $W^{1,1}(\Omega)$, and
associated with the norm
\begin{align*}
\| u\|_{W^{1,G}(\Omega)}= \|  u\|_{L^{G}(\Omega)}+\| \nabla u\|_{L^{G}(\Omega)}.
\end{align*}
We present some properties of spaces $L^{G}(\Omega)$ and $W^{1,G}(\Omega)$, and properties of functions in $L^{G}(\Omega)$ and $W^{1,G}(\Omega)$.
\begin{lemma}[\cite{MW}]\label{modular} There exists a constant $C=C({\delta_0},g_0)$ such that
\begin{align*}\| u\|_{L^{G}(\Omega)}\leq C \max \left\{\left(\int_{\Omega} G(|u|)\text{d}x \right)^{\frac{1}{1+{\delta_0}}}, \left(\int_{\Omega} G(|u|)\text{d}x\right)^{\frac{1}{1+g_0}}  \right\}. \end{align*}
\end{lemma}
\begin{lemma}[\cite{MW}] $L^{\tilde{G}}(\Omega)$ is the dual of $L^{G}(\Omega).$ Moreover, $L^{G}(\Omega)$ and $W^{1,G}(\Omega)$ are reflexive.
\end{lemma}
\begin{lemma}[\cite{MW}]\label{embedding}
$L^{G}(\Omega) \hookrightarrow L^{1+{\delta_0}}(\Omega)$ continuously.	
\end{lemma}
\begin{lemma}[\cite{MW}]\label{holder} For any $u \in L^{G}(\Omega)$ and any $v \in L^{\tilde{G}}(\Omega),$ there holds $\left|\int_{\Omega} uv \text{d}x \right|$ $ \leq 2 \| u\|_{L^{G}(\Omega)} \| v\|_{L^{\tilde{G}}(\Omega)}.$
\end{lemma}
\begin{lemma}[\cite{G}]\label{poincare} For any $u \in W^{1,G}_{0}(\Omega),$ which is the closure of $C^{\infty}_{0}(\Omega)$ in $W^{1,G}(\Omega),$ there holds $\int_{\Omega} G(|u|) \text{d}x \leq  \int_{\Omega} G(c|\nabla u|) \text{d}x$, where the constant $c$ is twice the diameter of $\Omega.$
\end{lemma}
\begin{lemma}[\cite{MW}]\label{continuity} Let $u\in L^\infty(\Omega)$ such that for some $\alpha \in (0,1)$ and $r_0>0$,
\begin{align*}
\int_{B_{r}}G(|\nabla u|)\text{d}x\leq & C_1r^{N+\alpha-1},\ \ 0<r\leq r_0,
\end{align*}
with $B_{r_0}\Subset \Omega$. Then $u\in C^\alpha(\Omega)$ and there exists a constant $C=C(C_1,\alpha,N,g_0,G(1))$ such that $[u]_{0,\alpha,\Omega}\leq C$.
\end{lemma}
\begin{lemma}\label{Lemma 2.131}
Let $v$ be a bounded weak solution of $\text{div} \ \frac{g(|\nabla v|)}{|\nabla v|}\nabla v=0$
 in $B_{R}$ {(see \eqref{weak solution} for a definition)}. For every $\lambda\in (0,N)$, there exists $C=C(\lambda,N,\delta,g_0,\|v\|_{L^{\infty}(B_{R})})>0$ such that
 \begin{align*}
\int_{B_{r}}G(|\nabla v|)\text{d}x\leq Cr^{\lambda},\ \ \ \forall 0<r\leq R.
\end{align*}
\end{lemma}
\begin{pf*}{Proof}
See \cite[(5.9), page 346]{L}, or \cite[Lemma 2.7]{MW}.\hfill $\blacksquare$
\end{pf*}
\begin{lemma}[{\cite{ZFZ}}]\label{Lemma +2.4}
 Let $u\in W^{1,G}(\Omega),
 B_{R}\subset\Omega$. If $v$ is a bounded weak solution of
\begin{align}
\text{div} \ \frac{g(|\nabla v|)}{|\nabla v|}\nabla v=0\ \  \text{ in}\  B_{R},\ \ v-u\in
W^{1,G}_{0}(B_{R}),\notag
\end{align}
then for any $\lambda\in (0,N)$, there exists $C=C(\lambda,N,{\delta_0},g_0,\|v\|_{L^{\infty}(B_{R})})>0$ such that
\begin{align*}
\int_{B_{R}}G(|\nabla u-\nabla v|)\text{d}x\leq  C\int_{B_{R}}(G(|\nabla u|)-G(|\nabla v|))\text{d}x+CR^{\frac{\lambda}{2}}\bigg(\int_{B_{R}}(G(|\nabla u|)-G(|\nabla v|))\text{d}x\bigg)^{\frac{1}{2}}.
\end{align*}
\begin{pf*}{Proof}
See the proof of Lemma 3.1 in \cite{ZFZ}.\hfill $\blacksquare$
\end{pf*}
\end{lemma}
 Let $(u)_{r}=\frac{1}{|B_{r}|}\int_{B_{r}}u\text{d}x$ be the
 average value of function $u$ on the ball $B_{r}$.
 \begin{lemma}[{\cite{ZFZ}}]\label{Lemma 4.1}\ \
 Let $u\in W^{1,G}(\Omega),
 B_{R}\subset\Omega$. If $v\in W^{1,G}(B_{R})$ is a weak solution of $\text{div} \ \frac{g(|\nabla v|)}{|\nabla v|}\nabla v=0$ in
$B_{R}$, then for some positive constant $0<\sigma<1$, there exists a positive
constant $C=C(N,\delta,g_{0})$ such that for each $0<r\leq R$, there
holds
\begin{align}
\int_{B_{r}} G(|\nabla u-(\nabla u)_{r}|)\text{d}x\leq
C\bigg(\frac{r}{R}\bigg)^{N+\sigma}\int_{B_{R}} G(|\nabla u-(\nabla
u)_{R}|)\text{d}x
+C\int_{B_{R}} G(|\nabla u-\nabla v|)\text{d}x.\notag
\end{align}
\end{lemma}
The following result is an iteration lemma, which will be used in the establishment of regularities of minimizers of $\mathcal {J}$.
 \begin{lemma}[{\cite{LdT}}]\label{Lemma 2.13}
\ \ Let $\overline{\phi}(s)$ be a
 non-negative and non-decreasing function. Suppose that
\begin{align}
\overline{\phi}(r)\leq
C_{1}\bigg(\bigg(\frac{r}{R}\bigg)^{\alpha}+\vartheta\bigg)\overline{\phi}(R)+C_{2}R^{\beta},\notag
\end{align}
for all $r\leq R\leq R_{0}$, with $C_{1},\alpha,\beta$ positive
constants and $C_{2},\vartheta$ non-negative constants. Then, for
any $\tau<\min\{\alpha,\beta\}$, there exists a constant
$\vartheta_{0}=\vartheta_{0}(C_{1},\alpha,\beta,\tau)$ such that if
$\vartheta<C_{1},\vartheta_{0}$, then for all $r\leq R\leq R_{0}$ we
have
\begin{align}
\overline{\phi}(r)\leq
C_{3}\bigg(\frac{r}{R}\bigg)^{\tau}\big(\overline{\phi}(R)+C_{2}R^{\tau}\big),\notag
\end{align}
where $C_{3}=C_{3}(C_{1},\tau-\min\{\alpha,\beta\})$ is a positive
constant. In turn,
\begin{align}
\overline{\phi}(r)\leq C_{4}r^{\tau},\notag
\end{align}
where $C_{4}=C_{4}(C_{2},C_{3},R_{0},\overline{\phi},\tau)$ is a
positive constant.
\end{lemma}

\section{Existence, $L^\infty-$boundedness and continuity of minimizers over the set $\mathcal {K}$}\label{Sec: existence and continuity}
\begin{theorem} \label{existence-boundness} Under the growth condition \eqref{hip-F} and $1+f_0 \leq (1+\delta_0)^{\star}$ where $(1+\delta_0)^{\star} := \frac{N(1+\delta_0)}{N-(1+\delta_0)}$ there exists a minimizer $u$ of the functional $\mathcal {J}(u)$ over the set $\mathcal {K} := \{u \in W^{1,G}(\Omega) ; u - \varphi \in W^{1,G}_{0}(\Omega)\}$, and there exists a constant $C_0>0$ depending only on ${\delta_0} , g_0, \lambda_{+}, {G(1)},\tilde{G}(1)$, {}{$\|h\|_{L^N(\Omega)}$, $\| \varphi\|_{L^{\infty}(\Omega)}$, and $\|u\|_{W^{1,G}(\Omega)} $}, such that
\begin{align*}\| u\|_{L^{\infty}(\Omega)} \leq C_0 , \end{align*}
for all minimizers $u$ of $\mathcal{J}(u)$ over the set $\mathcal {K}$.
\end{theorem}

\begin{pf*}{Proof}
First it will be proved that
\begin{equation}\label{defmin}
I_0 := \inf_{ v \in \mathcal{K}} \mathcal{J}(v) > -\infty.
\end{equation}
For $v \in \mathcal{K} $ we have from the definition of $\mathcal{J}$ that
\begin{equation}\label{defJ}
\mathcal{J}(v) \geq \int_{\Omega} G(|\nabla v|)  + q F(v^{+}) - |h| |v| \text{d}x
\end{equation}

Using H\"{o}lder's inequality, $(\tilde{G}_2)$, $W^{1,1}_0(\Omega) \hookrightarrow L^{1^{\star}}(\Omega),$ where $1^{\star}= \frac{N}{N-1},$ and Lemma \ref{modular}, for an arbitrary $\tau>0$, there is a constant $C_{\tau}$ such that
{}{
\begin{align}
\int_{\Omega} |h||v| \text{d}x & \leq \|h\|_{L^N(\Omega)}\|v\|_{L^{1^*}(\Omega)} \notag\\
& \leq \|h\|_{L^N(\Omega)}(\|v-\varphi\|_{L^{1^*}(\Omega)}+\|\varphi\|_{L^{1^*}(\Omega)})\notag\\
& \leq \|h\|_{L^N(\Omega)}(C\|\nabla v-\nabla\varphi\|_{L^{1}(\Omega)}+\|\varphi\|_{L^{1^*}(\Omega)})\notag\\
& \leq C\|h\|_{L^N(\Omega)}(\|\nabla v\|_{L^{1}(\Omega)} + \|\nabla\varphi\|_{L^{1}(\Omega)}+\|\varphi\|_{L^{1^*}(\Omega)})
\notag\\
& \leq C\|h\|_{L^N(\Omega)}(C_\tau\|1\|_{L^{\tilde{G}}(\Omega)}+\tau \|\nabla v\|_{L^{G}(\Omega)}+ \|\nabla\varphi\|_{L^{1}(\Omega)}+\|\varphi\|_{L^{1^*}(\Omega)})
\notag\\
& \leq K+K\tau\max \left\{\left(\int_{\Omega} G(|\nabla v|) \text{d}x\right)^{\frac{1}{1+\delta_0}}, \left(\int_{\Omega} G(|\nabla v|) \text{d}x\right)^{\frac{1}{1+g_0}}\right \}\notag\\
&\leq K + K \tau \left( 1 + \int_{\Omega} G(|\nabla v|) \text{d}x \right),\label{sobolev}
\end{align}
}
where $K>0$ is a constant that does not depend on $\tau.$

For $v \in \mathcal{K} $, we have from Lemma \ref{poincare} that
\begin{align}\label{poinc}
\int_{\Omega} G(|v-\varphi|) \text{d}x \leq \int_{\Omega} G (c |\nabla v - \nabla \varphi|) \text{d}x, c = 2 \text{diam}(\Omega).
\end{align}
Thus, using {$(G_3),(G_4), \eqref{poinc}$} and the nondecreasing monotonicity of $G$ we have
\begin{align}\label{aux-1}
\int_{\Omega} G(|v|) \text{d}x \leq &  \int_{\Omega} G(|v-\varphi| + |\varphi|) \text{d}x\notag\\
\leq &  C \left( \int_{\Omega} G(|v-\varphi|) \text{d}x + \int_{\Omega} G(|\varphi|) \text{d}x\right)\notag\\
\leq &  C \left( \int_{\Omega}G(|\nabla v - \nabla \varphi|) \text{d}x + \int_{\Omega} G(|\varphi|) \text{d}x\right)\notag\\
\leq& C \left( \int_{\Omega}G(|\nabla v|)\text{d}x + \int_{\Omega} G(|\nabla \varphi|)  \text{d}x + \int_{\Omega} G(|\varphi|) \text{d}x\right),
\end{align}
where $C$ is constant depending only on the diameter of $\Omega,$ and ${\delta_0}$ and $g_0 .$

The {}{hypothesis} $(f_2)$ implies that for given $\tau >0$, there exists a constant $K_\tau >0$ such that
\begin{align}\label{tau}
F(t) \leq K_{\tau} + \tau G(t), \ \ \ \  \forall t \geq 0.
\end{align}

By \eqref{tau}, we have
\begin{equation}\label{qF}
\left| \int_{\Omega}q F(v^+) \text{d}x\right| \leq \| q\|_{L^{\infty}(\Omega)}\int_{\Omega} {F(v^+)} \text{d}x \leq \| q\|_{L^{\infty}(\Omega)}K_{\tau} + \| q\|_{L^{\infty}(\Omega)}\tau \int_{\Omega} G(|v|) \text{d}x.
\end{equation}
Using \eqref{sobolev} ,\eqref{aux-1}, \eqref{qF} and considering a suitable choice of $\tau$ we obtain that $I_0 > -\infty.$

Now consider $v_j ( j \in \mathbb{N})$, a minimizing sequence, and $j_0 \in \mathbb{N}$, such that $\mathcal{J}(v_j) \leq I_0 + 1$ for all $j \geq j_0 .$
Arguing as in \eqref{sobolev} ,\eqref{aux-1} and \eqref{qF} we have for $\tau >0$ that
{}{\begin{align}
\int_{\Omega} G(|\nabla v_j|) \text{d}x =& \mathcal{J}(v_j) - \int_{\Omega}  q F(v^{+}) \text{d}x - \lambda_{+} \chi_{\{v>0\}}  - h v_j \text{d}x\notag\\
 \leq& \tau C \left( \int_{\Omega} G(|\nabla v_j|)\text{d} x+ G(|\nabla \varphi|) + G(|\varphi|)\text{d}x  \right) +K\label{tau-sequence}
\end{align}}
for all $ j \geq j_0,$ where {}{$C,K$ are constants} independent of $j \in \mathbb{N}.$ A suitable choice of $\tau>0$ implies that the sequence $\int_{\Omega} G(|\nabla v_j|) \text{d}x, j \in \mathbb{N}$ is bounded. The reasoning of  \eqref{aux-1} implies that the sequence  $\int_{\Omega} G(|v_j|) \text{d}x, j \in \mathbb{N}$ is bounded. Thus the sequence $(v_j), j \in \mathbb{N}$ is bounded in $W^{1,G}(\Omega).$

From $(G_4)$ and Lemma \ref{poincare}, we deduce that $v_j-\varphi$ is a bounded sequence in $W^{1,G}_{0}(\Omega).$ Since $W^{1,G}_{0}(\Omega)$ is reflexive, there exists $u \in W^{1,G}(\Omega)  $ with $u-\varphi\in W^{1,G}_0(\Omega)$ such that for a subsequence we have
\begin{align*}
 v_j \rightharpoonup u \ \text{in} \ W^{1,G}(\Omega).\end{align*}
Then by Lemma \ref{embedding}, we find that
\begin{align*}
v_j \rightharpoonup u \ \ \text{in} \ W^{1,1+{\delta_0}}(\Omega).
\end{align*}
Thus, up to a subsequence, we have that $v_j \rightarrow u$ a.e in $\Omega.$

Note that
\begin{align}\label{lim-inf}
\int_{\Omega}G(|\nabla u |) \text{d}x \le \liminf_{j \rightarrow \infty} \int_{\Omega}G(|\nabla v_j |) \text{d}x.
\end{align}
In fact, by the convexity of $G,$ it follows
\begin{align}\label{convexity}
\int_{\Omega} G(|\nabla v_j|) \text{d}x \geq \int_{\Omega} G(|\nabla u|)\text{d}x + \int_{\Omega} g(|\nabla u|)\frac{\nabla u}{|\nabla u|}(\nabla v_j - \nabla u) \text{d}x.
\end{align}
We have that
\begin{align*} \tilde{G}\left(g(|\nabla u|) \frac{\partial u }{\partial x_i} \frac{1}{|\nabla u|}\right) \leq \tilde{G}(g(|\nabla u|)) \leq C G(|\nabla u|),\end{align*}
which implies that $g(|\nabla u |) \frac{\nabla u}{|\nabla u|} \in {}{(L^{\tilde{G}}(\Omega))^N}.$ Thus, combining the fact that
$\nabla v_j \rightharpoonup \nabla u \ \text{in} \ {}{(L^{G}(\Omega))^N}$
with the inequality \eqref{convexity}, we get \eqref{lim-inf}. Since the sequence ${}{\{qF(v_j)\}}, j \in \mathbb{N}$, is bounded below, we have from Fatou's Lemma that
\begin{align}\label{Fatou}
\liminf_{j \rightarrow +\infty} \int_{\Omega} qF((v_j)^+) \text{d}x \geq  \int_{\Omega}q F(v^+) \text{d}x.
\end{align}
Note also that
\begin{align}\label{liminf}
\liminf_{j \rightarrow +\infty} \int_{\Omega}\lambda_{+}\chi_{\{v_j > 0\}} \text{d}x \geq \int_{\Omega}\lambda_{+}\chi_{\{u > 0\}} \text{d}x.
\end{align}
The continuous embedding $W^{1,G}(\Omega) \hookrightarrow L^{\frac{N}{N-1}}(\Omega)$ implies that $\{v_j\},j \in \mathbb{N}$, is a bounded sequence in $L^{\frac{N}{N-1}}(\Omega).$ Therefore
\begin{align}\label{lebesgue3}
\int_{\Omega}hv_j  \text{d}x \rightarrow  \int_{\Omega} h u \text{d}x.
\end{align}
Thus from \eqref{lim-inf}, \eqref{Fatou}, \eqref{liminf} and \eqref{lebesgue3} we deduce that $\mathcal{J}(u) \leq \liminf\limits_{j \rightarrow +\infty}\mathcal{J}(v_j),$ which implies that $u$ is a minimizer of $\mathcal{J}(u).$

Now we prove the boundedness of minimizers. Let $k_0 \in \mathbb{N}$ such that
\begin{equation}\label{defk_0}
k_0 \geq \| \varphi\|_{L^{\infty}(\partial \Omega)}.
\end{equation} For each $k \geq k_0$, define the function $u_k : \Omega \rightarrow \mathbb{R} $ by
\[
u_k =
\begin{cases}
k \cdot sgn(u) &\quad\text{if} \ |u|>k,\\
u  &\quad\text{if} \ |u| \leq k,
\end{cases}
\]
where $sgn(u) = 1$ if $u \geq 0 $, and $sgn(u)=-1$ if $u<0$. Let $A_k := \{|u| > k\}.$ For each $k \geq k_0$, we have
\begin{align*}u=u_k \ \text{in} \ {A_{k}}^c,  \ \ \ \text{and} \ \ \ \ u_k = k\cdot sgn(u) \  \text{in} \ A_{k}.\end{align*}
We have  $(|u|- k)^+\in W^{1,G}_{0}(\Omega) $ for all $k \geq k_0.$ Note that
\begin{align}\label{bound1}
\int_{A_{k}} G(|\nabla u|) \text{d}x =& \int_{\Omega} G(|\nabla u|) - G(|\nabla u_k|) \text{d}x  \notag\\
&\leq \int_{A_k} q(F(u_k^+) - F(u^+)) + \lambda_{+}(\chi_{\{u_k >0\}} - \chi_{\{u>0\}}) + h(u_k - u)\text{d}x\notag\\
&\leq \| q\|_{L^{\infty}(\Omega)}\int_{A_k} |F(u_k^+) - F(u^+)| \text{d}x +\int_{A_k} h(u_k - u) \text{d}x.
\end{align}
Note that
\begin{equation}\label{diferent}
\begin{aligned}
\int_{A_k} h (u_k - u) \text{d}x = & \int_{A_k \cap \{u \geq 0\}} h (k - u) \text{d}x + \int_{A_k \cap \{u<0\}} h (-k - u) \text{d}x \\
\leq & 2 \int_{A_k} |h| (|u| -k) \text{d}x.
\end{aligned}
\end{equation}
Using the continuous  embeddings $W^{1,G}(\Omega) \hookrightarrow W^{1,1}(\Omega),$ $W^{1,1}_{0}(\Omega) \hookrightarrow L^{1^{\star}}(\Omega),$ $(\tilde{G}_2)$ and the H\"{o}lder's inequality we get
\begin{equation}\label{dif-2}
\begin{aligned}
\int_{A_k} |h| (|u| -k) \text{d}x \leq & C \left(\int_{A_k} |h|^N \text{d}x\right)^{\frac{1}{N}} \left( \int_{A_k} (|u| - k)^{\frac{N}{N-1}} \text{d}x\right)^{\frac{N-1}{N}}\\
\leq & C \left(\int_{A_k} |h|^N \text{d}x\right)^{\frac{1}{N}}  \left( \int_{\Omega} ((|u|-k)^{+})^{\frac{N}{N-1}} \text{d}x \right)^{\frac{N-1}{N}}\\
\leq & C \int_{\Omega} |\nabla (|u|-k)^{+}| \text{d}x \\
=& C \int_{A_k} |\nabla  u| \text{d}x \\
\leq & C \left( \varepsilon \int_{A_k} G(|\nabla u|) \text{d}x +|A_k|\right),
\end{aligned}
\end{equation}
where $C$ is a constant that does not depend on $\varepsilon$ and $k.$

{}{Note that by $(f_{1})$, there exists $\alpha,\beta$ satisfying
\begin{align*}
1<1+\alpha \leq (1+\delta_0)^{\star},\ \text{and}\  0<\beta<1,
\end{align*}
such that
\begin{align}
f(t)\leq f(1)\frac{1+f_0}{1+\theta_0}(1+t^\alpha+t^{-\beta}),\forall t> 0.\label{+1114}
\end{align}}
The Mean Value Theorem imply for  some $\theta \in (0,1),$ which depends on $x \in \Omega $ and $k \in \mathbb{N}$, that
\begin{align}\label{F-1}
\int_{A_k} |F(u_k^+) - F(u^+)| \text{d}x=&\int_{A_k\cap\{u\geq0\}} (F(u) - F(k) ) \text{d}x \notag\\
=&\int_{A_k\cap\{u\geq0\}} (F(|u|) - F(k)) \text{d}x \notag\\ \leq & \int_{A_k\cap\{u\geq 0\}} f( (1-\theta)k + \theta |u|)(|u|-k)  \text{d}x\notag\\
\leq &C \int_{A_k\cap\{u\geq 0\}}(1 + |(1-\theta)k + \theta |u||^{-\beta} + |(1-\theta)k + \theta |u||^{\alpha})(|u|-k) \text{d}x \notag \\
\leq & C \int_{A_k \cap \{u \geq 0\}}(1 + |(1-\theta)k + \theta k|^{-\beta} + |(1-\theta)k + \theta |u||^{\alpha})(|u|-k) \text{d}x\notag \\
\leq& C \int_{A_k \cap \{u \geq 0\}} (1 + k+ |u|)^{\alpha}(|u|-k) \text{d}x \notag \\
\leq& C \int_{A_k \cap \{u \geq 0\}} (1 +  |u|)^{\alpha}(|u|-k) \text{d}x \notag \\
\leq&  C \left( \int_{A_k} ||u|-k| \text{d}x +  \int_{A_k} |u|^{\alpha}||u|-k| \text{d}x \right)
\end{align}
The embeddings $W^{1,G}(\Omega) \hookrightarrow W^{1,1+\delta_0}(\Omega) \hookrightarrow L^{(1+\delta_0)^{\star}}(\Omega)$ implies that $|u|^{\alpha} \in L^{\frac{1+\alpha}{\alpha}}(\Omega)$ and $||u|-k| \in L^{1+\alpha}(\Omega).$
The Young's inequality implies that
\begin{align}\label{F-2}
\int_{A_k}|u|^{\alpha} ||u| -k| \text{d}x \leq & C \left( \int_{A_k} |u|^{1+\alpha} \text{d}x + \int_{A_k} ||u| - k|^{1+\alpha} \text{d}x\right) \notag \\
\leq & C \left( \int_{A_k} ||u|-k|^{1+\alpha} + k^{1+\alpha} \text{d}x + \int_{A_k} ||u| - k|^{1+\alpha} \text{d}x\right) \notag \\
\leq & C \left( \int_{A_k} ||u| - k|^{1+\alpha} \text{d}x + k^{1+\alpha} |A_k| \right)
\end{align}
Thus from \eqref{F-1} and \eqref{F-2} we have
\begin{equation}\label{Flast}
\| q\|_{L^{\infty}(\Omega)} \int_{A_k} |F(u_k^{+}) - F(u^{+})| \text{d}x \leq C \left( \int_{A_k} ||u| - k|^{1+\alpha} \text{d}x + k^{1+\alpha} |A_k| \right)
\end{equation}

Considering  a suitable choice of $\varepsilon$ in \eqref{dif-2}, using    \eqref{diferent}\ and \eqref{Flast}  we obtain that
\begin{equation}\label{Ffinal}
\int_{A_k} G(|\nabla u| ) \text{d}x \leq C \left( \int_{A_k} ||u|-k|^{1+\alpha} \text{d}x + k^{1+\alpha} |A_k|\right).
\end{equation}
From $(G_2)$ we have
\begin{align}\label{Ffinal2}
\int_{A_k} G(|\nabla u|) \text{d}x  =&  \int_{A_k \cap \{|\nabla u| \leq 1|\}} G(|\nabla u|)\text{d}x + \int_{A_k \cap \{|\nabla u| >1\} } G(|\nabla u|) \text{d}x \notag \\
\geq & C \left( \int_{A_k \cap \{|\nabla u| \leq 1|\}} |\nabla u|^{1+g_0}\text{d}x + \int_{A_k \cap \{|\nabla u| >1\} } |\nabla u|^{1+\delta_0} \text{d}x \right).
\end{align}
We obtain from \eqref{Ffinal} and \eqref{Ffinal2}
\begin{equation}
\int_{A_k} |\nabla u|^{1+\delta_0} \text{d}x \leq  C \left(\int_{A_k} ||u|-k|^{1+\alpha} \text{d}x  + (k^{1+\alpha}+1)|A_k| \right).
\end{equation}
Since $\nabla |u| = (\nabla u) sgn(u)$ we have $|\nabla u| = |\nabla |u||$. Therefore
\begin{equation}\label{fusco}
\int_{A_k} |\nabla |u||^{1+\delta_0} \text{d}x \leq C \left(\int_{A_k} ||u|-k|^{1+\alpha} \text{d}x  + (k^{1+\alpha} +1)|A_k| \right)
\end{equation}

Define
\begin{equation}\label{defK}
K_n := \displaystyle\frac{K}{2} \left( 1 - \displaystyle\frac{1}{2^{n+1}}\right), n=0,1,2,...
\end{equation}
and suppose that $K_n , K \geq k_0 , n=0,1,2,..., $ where $k_0$ is given in \eqref{defk_0}. Consider
\begin{align*}J_n:= \int_{A_{K_n}}  ((|u|-K_n)^{+})^{(1+\delta_0)^{\star}} \text{d}x, n=0,1,2,... .\end{align*}
We claim that $J_{n+1} \leq C D^{n} {J_n}^{1+\xi}, n=0,1,2,... $ with $C,\xi >0$ and $D>1 $ with $C,D$ and $\xi$ not depending on $n.$
The continuous embedding $W^{1+\delta_0}_0(\Omega) \hookrightarrow L^{(1+\delta_0)^{\star}}(\Omega)$ combined with the fact that $(|u|- K_{n})^{+} \in W^{1,1+\delta_0}_0 (\Omega), n=0,1,...$ implies that
\begin{align}\label{fusco2}
J_{n+1} \leq & \int_{A_{K_{n+1}}} ((|u|- K_{n+1})^{+})^{(1+\delta_0)^{\star}} \text{d}x \notag \\
\leq & \int_{\Omega}  ((|u|- K_{n+1})^{+})^{(1+\delta_0)^{\star}} \text{d}x \notag \\
\leq & C\left( \int_{\Omega} |\nabla (|u|-K_{n+1})^{+}|^{1+\delta_0} \text{d}x\right)^{\frac{(1+\delta_0)^{\star}}{1+\delta_0}} \notag \\
=& C\left( \int_{A_{K_{n+1}}} |\nabla |u||^{1+\delta_0} \text{d}x\right)^{\frac{(1+\delta_0)^{\star}}{1+\delta_0}}
\end{align}
From \eqref{fusco} and \eqref{fusco2} we get
\begin{equation}\label{fusco3}
J^{\frac{1+\delta_0}{(1+\delta_0)^{\star}}}_{n+1} \leq C \left( \int_{A_{K_{n+1}}} ||u|-K_{n+1}|^{(1+\delta_0)^{\star}} \text{d}x + (K^{(1+\delta_0)^{\star}}_{n+1} +1)|A_{k_{n+1}}|\right)
\end{equation}
We have $K_{n+1} - K_n = \displaystyle\frac{K}{2^{n+3}},$ therefore
\begin{align*}
\left( \frac{K}{2^{n+3}}\right)^{(1+\delta_0)^{\star}} |A_{K_{n+1}}| = & (K_{n+1} - K_n)^{(1+\delta_0)^{\star}} |A_{K_{n+1}}| \notag \\
=& \int_{A_{K_{n+1}}} |K_{n+1} - K_n|^{(1+\delta_0)^{\star}}\text{d}x \notag \\
\leq & \int_{A_{K_{n+1}}} ((|u| - K_n)^{+})^{(1+\delta_0)^{\star}} \text{d}x \notag \\
\leq & J_n,
\end{align*}
Therefore
\begin{equation}\label{fusco4}
|A_{K_{n+1}}| \leq \left(\frac{2^{n+3}}{K}\right)^{(1+\delta_0)^{\star}} J_n
\end{equation}
From \eqref{fusco4} and using that $K \geq 1$ we get
\begin{align}\label{fusco5}
\int_{A_{K_{n+1}}} ||u|-K_{n+1}|^{(1+\delta_0)^{\star}} \text{d}x \leq& \int_{A_{K_{n+1}}} ||u|-K_{n}|^{(1+\delta_0)^{\star}} \text{d}x +  \int_{A_{K_{n+1}}} |K_n - K_{n+1}|^{(1+\delta_0)^{\star}} \text{d}x \notag \\
\leq &\int_{A_{K_{n}}} ||u|-K_{n}|^{(1+\delta_0)^{\star}} \text{d}x +  |K_n - K_{n+1}|^{(1+\delta_0)^{\star}} |A_{K_{n+1}}| \notag \\
\leq & J_n + |K_{n+1} - K_n|\left(\frac{2^{n+3}}{K}\right)^{(1+\delta_0)^{\star}} J_n  \notag \\
=& J_n + \left( \frac{K}{2^{n+3}}\right) \left(\frac{2^{n+3}}{K}\right)^{(1+\delta_0)^{\star}} J_n  \notag \\
\leq & J_n + 2^{(1+\delta_0)^{\star}(n+3)} J_n
\end{align}
We also have
\begin{align}\label{fusco6}
(K^{(1+\delta_0)^{\star}}_{n+1} +1)|A_{k_{n+1}} | = & \left( \left(\frac{K}{2} \left( 1 - \frac{1}{2^{n+2}}\right) \right)^{(1+\delta_0)^{\star}} +1\right)\left(\frac{2^{n+3}}{K}\right)^{(1+\delta_0)^{\star}} J_n  \notag \\
\leq & C 2^{(1+\delta_0)^{\star}(n+3)}J_n,
\end{align}
where $C>0$ is a constant that does not depend on $n.$
From \eqref{fusco5} and \eqref{fusco6} we have $J_{n+1} \leq C D^n J^{\frac{(1+\delta_0)^{\star}}{(1+\delta_0)}}_n $ where  $C>0$ and $D>1$ are constants that does not depend on $n.$ Therefore $J_{n+1} \leq C D^n J^{1+\xi}_n$ with $C,\xi>0$ and $D>1$ are constants that does not depend on $n.$
Since
\begin{align*}J_0 : = \int_{A_0} ((|u | - K_0)^{+})^{(1+\delta_0)^{\star}}\text{d}x \leq  \int_{\Omega} \left(\left(|u|-\frac{K}{2}\right)^{+}\right)^{(1+\delta_0)^{\star}}\text{d}x,\end{align*}
it is possible to choose $K$ large enough, that depends on {}{$\|u\|_{L^{(1+\delta_0)^{\star}}(\Omega)}$},  such that $J_0 \leq C^{-\frac{1}{\xi}} D^{-\frac{1}{\xi^2}}$ . By Lemma \ref{lady-aux} we have $J_n \rightarrow 0$ as $n \rightarrow +\infty.$ On other hand we have
\begin{align*}\lim_{n \rightarrow +\infty} \int_{A_{K_n}} ((|u|-K_n)^{+})^{(1+\delta_0)^{\star}} \text{d}x = \int_{A_{\frac{K}{2}}} \left(\left(|u|-\frac{K}{2}\right)^{+}\right)^{(1+\delta_0)^{\star}} \text{d}x,\end{align*}
therefore $|u| \leq K$ a.e in $\Omega,$  where $K$ is a constant that depends on {}{$\|u\|_{L^{(1+\delta_0)^{\star}}(\Omega)}$, and therefore, depends on $\|u\|_{W^{1,G}(\Omega)}$.}
 $\blacksquare$

\end{pf*}

{}{In the next result we prove that the minimizers of $\mathcal{J}$ are uniformly bounded under some conditions.}

\begin{theorem}[Uniform $L^{\infty}$ estimate] \label{unif-boundness}Under the growth condition \eqref{hip-F} and $1+f_0 < (1+\delta_0)^{\star}$ where $(1+\delta_0)^{\star} := \frac{N(1+\delta_0)}{N-(1+\delta_0)}$, there exists a constant $C_0 >0$ depending only on ${\delta_0} , g_0, \lambda_{+}, {G(1)},\tilde{G}(1)$, {}{$\|h\|_{L^N(\Omega)}$,$\| \varphi\|_{L^{\infty}(\Omega)}$, and $\|\varphi\|_{W^{1,G}(\Omega)}$}, such that
{}{\begin{align*}
\| u\|_{W^{1,G}(\Omega)} \leq C_0,
 \end{align*}}
 for all minimizers $u$ of $\mathcal{J}(u)$ over the set $\mathcal {K}$.
 Furthermore, there exists a constant $C_1 >0$ depending only on ${\delta_0} , g_0, \lambda_{+}, {G(1)},\tilde{G}(1)$, {}{$\|h\|_{L^N(\Omega)}$,$\| \varphi\|_{L^{\infty}(\Omega)}$, and $\|u\|_{W^{1,G}(\Omega)}$}, such that
\begin{align*}
\| u\|_{L^{\infty}(\Omega)} \leq C_1 ,
\end{align*}
for all minimizers $u$ of $\mathcal{J}(u)$ over the set $\mathcal {K}$.
\end{theorem}

\begin{pf*}{Proof}
Consider $\{v_j\}, j \in \mathbb{N}$ a sequence of minimizers with $\| v_j\|_{L^{\infty}(\Omega)} \rightarrow +\infty$  as $j \rightarrow +\infty.$ Using the fact that $\mathcal{J}(v_j) = I_0,$ where $I_0$ is given by \eqref{defmin}, and the reasoning that provides \eqref{tau-sequence} we get for an arbitrary $\tau >0$ that
\begin{equation}\label{unif-tau}
\int_{\Omega} G(|\nabla v_j|) \text{d}x \leq \tau C \left( \int_{\Omega} G(|\nabla v_j|) + G(|\nabla \varphi|) +  G(|\varphi|) \text{d}x \right)+C,
\end{equation}
where $C$  is a constant that does not depend on $\tau$ and $j \in \mathbb{N}.$ A suitable choice of $\tau $ in \eqref{unif-tau} and implies that the sequence $\int_{\Omega} G(|\nabla v_j|) \text{d}x$ is bounded. Thus arguing as in \eqref{aux-1} we obtain that the sequence $\{v_j\}, j \in \mathbb{N}$  is bounded in $W^{1,G}(\Omega).$

Since {}{$1+f_0 < (1+\delta_0)^{\star}$, we have $1+ \alpha < (1+\delta_0)^{\star}$ for some $\alpha>0$ in \eqref{+1114}. Then there is $\gamma\in (0, 1+\delta_0)$} such that $1+\alpha < \gamma^{\star} < (1+\delta_0)^{\star}$ where $\gamma^{\star}:= \frac{N\gamma}{N - \gamma }.$ By Lemma \ref{embedding}, we have the continuous embedding $W^{1,G}(\Omega) \hookrightarrow W^{1,1+\delta_0}(\Omega).$ Since the embedding $W^{1,G}(\Omega) \hookrightarrow L^{\gamma^{\star}}(\Omega)$ is compact it follows that there is $u \in W^{1,G}(\Omega)$ with $u-\varphi \in W^{1,G}_{0}(\Omega)$ such that $v_j \rightarrow u$ in $L^{\gamma^{\star}}(\Omega).$
For $k >k_0$ where $k_0$ is given by \eqref{defk_0} define the functions $u_{jk}: \Omega \rightarrow \mathbb{R}$ given by
\begin{equation}\label{defk_0}
k_0 \geq \| \varphi\|_{L^{\infty}(\partial \Omega)}.
\end{equation} For each $k \geq k_0$ and $j \in \mathbb{N}$ consider the functions $u_{jk} : \Omega \rightarrow \mathbb{R} $ given by
\[
u_{jk} =
\begin{cases}
k \cdot sgn(u_j) &\quad\text{if} \ |u_j|>k,\\
u_j  &\quad\text{if} \ |u_j| \leq k,
\end{cases}
\]
We have $(|u_{jk} - k|)^{+} \in W^{1,G}_0(\Omega)$ for all $k \geq k_0$ and $j \in \mathbb{N}.$ For an arbitrary $K>0$ define the sequence
\begin{align*}K_n := \frac{K}{2}\left( 1 - \frac{1}{2^{n+1}}\right), n=0,1,2,... .\end{align*}
Arguing as in the proof of Theorem \ref{existence-boundness} we obtain that
\begin{align*}\int_{A_{jk}} |\nabla |u||^{1+\delta_0} \text{d}x \leq C \left( \int_{A_{jk}} ||u_{jk}| - k|^{1+\alpha} \text{d}x + (k^{1+\alpha} +1)|A_{jk}|\right), \end{align*}
where $C>0$ is a constant that does not depend on $j$ and $k$ and $A_{jk}:= \{|u_j| > k\}.$
For each $j \in \mathbb{N}$ and  supposing that $k_n \geq K_0, n=0,1,2,...$ consider  the quantities
\begin{align*} J_{j,n+1}:=\int_{A_{j K_n}} ((|u_j| - k_n)^{+})^{\gamma^{\star}} \text{d}x.\end{align*}
Arguing as in \eqref{fusco2}, \eqref{fusco3}, \eqref{fusco4}, \eqref{fusco5} and \eqref{fusco6} we get for all $j$ and $n$ that
\begin{equation}\label{uniform-inequality}
J_{j,n+1} \leq CD^{n} J^{1+\xi}_{j,n},
\end{equation}
where $1+\xi = \frac{\gamma^{\star}}{\gamma}$ and $C$ and $D$ are constants that does not depend on $n$ and $j.$ Note that
\begin{align}\label{v0}
J_{j,0} &= \int_{A_{j k_0}} ((|v_j| - k_0)^{+})^{\gamma^{\star}} \text{d}x \notag \\
&\leq \int_{\Omega} \left(\left(|v_j| - \frac{K}{2}\right)^{+}\right)\text{d}x
\end{align}
Let $\varepsilon>0$ arbitrary. Since $v_j \rightarrow u$ in $L^{\gamma^{\star}}(\Omega)$ we have from \eqref{v0} that
\begin{align*}
J_{j,0} \leq \int_{\Omega} \left( \left(|u| - \frac{K}{2}\right)^{+}\right)^{\gamma^{\star}} + \frac{\varepsilon}{2}
\end{align*}
for all $j \geq j_0 ,$ where $j_0 \in \mathbb{N}$ is a number that depends only on $\varepsilon.$ Let $K>0$ large enough depending only on $u$ such that
\begin{align*}
\int_{\Omega} \left(\left( |u| - \frac{K}{2}\right)^{+}\right)^{\gamma^{\star}}\text{d}x < \frac{\varepsilon}{2},
\end{align*}
which implies that $J_{j,0} \leq \varepsilon$ for all $j \in \mathbb{N}$ with $j \geq j_0 .$ A suitable choice of $\varepsilon >0$ implies that  $J_{j,0} \leq C^{-\frac{1}{\xi}} D^{-\frac{1}{\xi^2}}$ for all $j \in \mathbb{N}$ with $j \geq j_0$ where $C,D$ and $\xi$ are that same constants of \eqref{uniform-inequality}. Thus by Lemma \ref{lady-aux} we have $J_{j,n} \rightarrow 0$ as $n \rightarrow +\infty$ for all $j \in \mathbb{N}$ with $ j \geq j_0.$ Therefore
\begin{align*}\int_{A_j\frac{K}{2}} \left(\left( |v_j| - \frac{K}{2}\right)^{+}\right)^{\gamma^{\star}}=0, \end{align*}
for all $j \in \mathbb{N}$ with $j \geq j_0.$ Thus $|v_j| \leq \frac{K}{2}$ a.e in $\Omega$ for all $j \in \mathbb{N}$ with $j \geq j_0 ,$ which proves the result.  $\blacksquare$
\end{pf*}
\begin{remark}
	Note that the proof of existence of minimizers holds in the case where $q$ is nonnegative and the function $F$ is only a continuous nonnegative function defined in $[0,+\infty).$
	We also point out that the $L^{\infty}$ estimates from Theorems \ref{existence-boundness}, and  \ref{unif-boundness}, hold if we replace condition $\eqref{hip-F}$ by the inequality
	\begin{align*}f(t) \leq C(1+ t^{-\beta} + t^{\alpha}), t>0 \end{align*}
	{}{where $\alpha $ and $\beta $ are constants with  $1<1+\alpha \leq (1+\delta_0)^{\star} $ and $0<\beta<1$, and $1<1+\alpha < (1+\delta_0)^{\star}$ and $0<\beta<1$, respectively, with $F(t) = \displaystyle\int_{0}^{t}f(s) \text{d}s,$ $t>0.$}
\end{remark}

\begin{corollary}[Local $C^{0,1-\frac{N}{m}}-$continuity]  \label{Holder continuity}
{}{Assume further that $h\in L^m(\Omega)$ with $m>N$. Let $u$ be a minimizer of $\mathcal{J}(u)$ over the set $\mathcal {K}$. Then $u \in C_{loc}^{0,1-\frac{N}{m}} (\Omega)$.}
\end{corollary}
\begin{pf*}{Proof}For any ball $B_R\Subset \Omega$, let $v$ be a weak solution of the following equation
\begin{align}
\left\{
  \begin{array}{ll}\nonumber
   \text{div}\ \frac{g(|\nabla v|)}{|\nabla v|}\nabla v=&0\ \ \ \ \text{in }\ B_{R}, \\
    \ \ \ \ \ \ \ \ \ \ \ \ \ \ \ \  \ \  \ \ \ v=&u\ \ \  \text{on}\ \partial B_{R}.
  \end{array}\notag
\right.
\end{align}
By the minimality of $u$, we have
\begin{align}\label{+15'}
\int_{B_{R}}G(|\nabla u|)\text{d}x-\int_{B_{R}}G(|\nabla v|)\text{d}x\leq&
\int_{B_{R}}q(F(v^+)-F(u^+))\text{d}x+\int_{B_{R}}h(v^+-u^+)\text{d}x
+\lambda_+\int_{B_{R}}(\chi_{\{v>0\}}
-\chi_{\{u>0\}})\text{d}x\notag\\
\leq & \|q\|_{L^\infty(B_{R})}\cdot\int_{B_{R}}(|F(v^+)|+|F(u^+)|)\text{d}x+C|B_R|^{\frac{m-1}{m}}\|h\|_{L^m(B_R)}+CR^N
\notag\\
\leq & \|q\|_{L^\infty(\Omega)}\cdot\int_{B_{R}}(F(\|v\|_{L^\infty(B_{R})})+F(\|u\|_{L^\infty(B_{R})}))\text{d}x+CR^{N-\frac{N}{m}}
\notag\\
\leq &  CR^{N-\frac{N}{m}},
\end{align}
where we used the increasing property of $F$, and the fact that $\|v\|_{L^\infty(B_R)}\leq \|v\|_{L^\infty(\partial B_R)}=\|u\|_{L^\infty(\partial B_R)}\leq \|u\|_{L^\infty(B_R)}\leq \|u\|_{L^\infty(\Omega)}\leq C$, which is guaranteed by the maximum principle.
By \eqref{+15'} and Lemma \ref{Lemma 2.131}, for any $\lambda \in (0,N)$, there holds
\begin{align*}
\int_{B_{R}}G(|\nabla u|)\text{d}x\leq & CR^N+\int_{B_{R}}G(|\nabla v|)\text{d}x\leq CR^{N-\frac{N}{m}}+CR^{\lambda}\leq CR^{N-\frac{N}{m}},
\end{align*}
where we let $N-\frac{N}{m}<\lambda<N$.
We conclude the desired result by Lemma \ref{continuity}.
\hfill $\blacksquare$
 \end{pf*}
 \begin{corollary}\label{Corollary 17}
 Assume further that $h\in L^\infty(\Omega)$ and $h\leq 0$ a.e. in $\Omega$, and $ \varphi \geq 0$ on $\partial \Omega$ in the sense of trace. Then every minimizer of $\mathcal{J}(u)$ over the set $\mathcal {K}$ is non-negative in $\Omega$.
\end{corollary}
\begin{pf*}{Proof}Let $\xi=\min\{u,0\}\leq 0$ and $\varepsilon \in (0,1)$. By the minimality of $u$, it follows
\begin{align*}\int_{\Omega} \bigg(G(|\nabla (u- \varepsilon \xi)|) -  G(|\nabla u|) + q\big(F((u-\varepsilon \xi)^{+})  - F(u^{+})\big)-h \varepsilon \xi+\lambda_{+} (\chi_{\{u - \varepsilon \xi >0\}}- \chi_{\{u >0\}})\bigg)\text{d}x\geq 0.
\end{align*}
Note that
\begin{align*}\int_{\Omega} q\big(F((u-\varepsilon \xi)^{+})  - F(u^{+})\big)\text{d}x=&\int_{\{u>0\}} q\big(F((u-\varepsilon \xi)^{+})  - F(u^{+})\big)\text{d}x+\int_{\{u\leq 0\}} q\big(F((u-\varepsilon \xi)^{+})  - F(u^{+})\big)\text{d}x\notag\\
=&\int_{\{u>0\}} q\big(F(u^{+})  - F(u^{+})\big)\text{d}x+\int_{\{u\leq 0\}} q\big(F((1-\varepsilon) u^{+})  - F(u^{+})\big)\text{d}x\notag\\
=& 0,
\end{align*}
and
\begin{align*}
\int_{\Omega} \lambda_{+} (\chi_{\{u - \varepsilon \xi >0\}}- \chi_{\{u >0\}})\text{d}x=&\int_{\{u>0\}} \lambda_{+} (\chi_{\{u  >0\}}- \chi_{\{u >0\}})\text{d}x+\int_{\{u\leq 0\}} \lambda_{+} (\chi_{\{(1 - \varepsilon)u >0\}}- \chi_{\{u >0\}})\text{d}x\notag\\
=&0.
\end{align*}
Therefore
\begin{align}\label{subsol'}
0 \leq & \frac{1}{\varepsilon}\int_{\Omega} \bigg(G(|\nabla (u- \varepsilon \xi)|) -  G(|\nabla u|) + q(F((u-\varepsilon \xi)^{+})  - F(u^{+}))-h \varepsilon \xi+\lambda_{+} (\chi_{\{u - \varepsilon \xi >0\}}- \chi_{\{u >0\}})\bigg)\text{d}x\\
\leq &-\int_{\Omega} g(|\nabla u -\varepsilon \nabla\xi|)\frac{\nabla u -\varepsilon \nabla\xi}{|\nabla u -\varepsilon \nabla\xi|}\nabla \xi\text{d}x,
\end{align}
where we used the convexity of $G$ and $h\xi \geq 0$ a.e. in $\Omega$.
Letting $\varepsilon \rightarrow 0^{+}$, we get
$
\int_{\Omega}\frac{g(|\nabla u|)}{|\nabla u|} \nabla u \nabla \xi \text{d}x  \geq 0.
$
It follows
\begin{align*}
\int_{\{u<0\}}g(|\nabla u|)|\nabla u| \text{d}x  \leq 0,
\end{align*}
which implies $u\equiv C $ or $u\geq 0$ a.e. in $\Omega$. By the fact that $u=\varphi\geq 0$ on $\partial \Omega$, and the continuity of $u$, we conclude that $u\geq 0$ in $\Omega$.
\hfill $\blacksquare$
\end{pf*}
\begin{remark}\label{Remark 4}
 Without restrictions on sign of $h$, it is easy to see by checking the proof of Theorem \ref{existence-boundness} and Corollary \ref{Holder continuity} that there exists a non-negative minimizer, which is also bounded and $C_{loc}^{0,1-\frac{N}{m}} -$contunous {}{under further assumption that $h\in L^m(\Omega)(m>N)$}, of the functional in \eqref{main equation} over the set $\widetilde{\mathcal {K}}=\{v \in W^{1,G}(\Omega): v-\varphi\in W^{1,G}_0(\Omega),v\geq 0\ a.e. \text{in}\ \Omega\}$ provided a non-negative $\varphi$.
\end{remark}
 \section{Local $C^{1,\alpha}-$ and Log-Lipschitz regularities of minimizers over the set $\mathcal {K}$}\label{Sec: 4}
In this section, we establish local $C^{1,\alpha}-$ and Log-Lipschitz continuities for minimizers of $\mathcal {J}(u)$.
We assume further that
\begin{align}
{}{1+f_0 < \displaystyle\frac{N(1+\delta_0)}{N-(1+\delta_0)}},f(t)=F'(t) \ \text{is\ monotone \ in} \ t>0.\label{1115+1}
\end{align}
 \begin{theorem}[Local $C^{1,\alpha}-$regularity of minimizers for $\lambda_+=0$]\label{Theorem +1.1}
{}{Assume that $ h \in L^{m}(\Omega)$ with $m>N$.} Let $u$ be a minimizer of $ \mathcal {J}(u)$ over the set $\mathcal {K}$ with $\lambda_+=0$. Then $u\in C^{1,\alpha}_{loc}(\Omega)$
for some $\alpha\in(0,1)$. More precisely, for any
$\Omega'\Subset\Omega$, there exists a constant $C>0$,
depending only on
$n,\theta_0,f_0,\delta_0,g_{0},G(1),\frac{1}{G(1)}$, {}{$\|h\|_{L^{m}(\Omega)}$}, $\|q\|_{L^{\infty}(\Omega)},\|\varphi
\|_{L^{\infty}(\Omega)}$, $\|\varphi
\|_{W^{1,G}(\Omega)}$ and $\Omega'$,
such that
\begin{align}
\|u\|_{C^{1,\alpha}(\Omega')}\leq C.\notag
\end{align}
  \end{theorem}
\begin{theorem}[Local Log-Lipschitz regularity of minimizers for $\lambda_+\geq 0$]\label{Theorem Log-Lip}
{}{Assume that $ h \in L^{N}(\Omega)$.} Let $u$ be a minimizer of $ \mathcal {J}(u)$ over the set $\mathcal {K}$  with $\lambda_+\geq 0$. Then $u$ is locally Log-Lipschitz continuous. More precisely, for any
$\Omega'\Subset\Omega$, there exists a constant $C>0$,
depending only on $n,\theta_0,f_0,\delta_0,g_{0},G(1),\frac{1}{G(1)}$, $\|h\|_{L^{N}(\Omega)}$, $\|q\|_{L^{\infty}(\Omega)},\|\varphi
\|_{L^{\infty}(\Omega)}$, $\|\varphi
\|_{W^{1,G}(\Omega)}$ and $\Omega'$, such that
\begin{align}
|u_{0}(x)-u_{0}(y)|\leq C|x-y||\log|x-y||,\notag
\end{align}
for any $x,y\in \Omega'$. Therefore,
$u_{0}\in C^{0,\tau}_{loc}(\Omega)$ for any $\tau<1$.
  \end{theorem}

 \begin{pf*}{Proof of Theorem \ref{Theorem +1.1}}\ Let $B_{R}=B_{R}(x_{0})$ for some $R\leq R_{0}\leq 1$,
  where $R_{0}$ will be chosen later. Without loss of generality, assume that $B_r\Subset B_R\Subset\Omega$, and $B_r$ and $B_R$ have the same centre. Let $v$ be a $G-$harmonic function
   in $B_{R}$ that agrees
 with $u$ on the boundary, i.e.,
\begin{align}
\text{div}\ \frac{g(|\nabla v|)}{|\nabla v|}\nabla v=0\ \text{in}\
B_{R}\ \ \text{and}\ \ v-u\in W^{1,G}_{0}(B_{R}).\notag
\end{align}
By Lemma \ref{Lemma 4.1} and Lemma \ref{Lemma +2.4},
we have
\begin{align}
\int_{B_{r}} G(|\nabla u-(\nabla u)_{r}|)\text{d}x\leq &
C\bigg(\frac{r}{R}\bigg)^{N+\sigma}\int_{B_{R}} G(|\nabla u-(\nabla
u)_{R}|)\text{d}x+C\int_{B_{R}} G(|\nabla u-\nabla v|)\text{d}x\notag\\
\leq & C\bigg(\frac{r}{R}\bigg)^{N+\sigma}\int_{B_{R}} G(|\nabla u-(\nabla
u)_{R}|)\text{d}x+C\int_{B_{R}}(G(|\nabla u|)-G(|\nabla v|))\text{d}x\notag\\
&+CR^{\frac{\lambda}{2}}\bigg(\int_{B_{R}}(G(|\nabla u|)-G(|\nabla v|))\text{d}x\bigg)^{\frac{1}{2}},\label{4.11}
\end{align}
where $\lambda $ is an arbitrary constant in $(0,n)$.

The minimality of $u$ and the increasing monotonicity of $F$ imply that
\begin{align}\label{5281}
\int_{B_{R}} (G(|\nabla u|)-G(|\nabla v|))\text{d}x\leq & \int_{B_{R}}
\bigg(q(F(v^+)-F(u^+))+h(v-u)\bigg)\text{d}x.\notag\\
 \leq& \|q\|_{L^{\infty}(B_R)}\int_{B_{R}}
|F(v^+)-F(u^+)|\text{d}x+{}{\int_{B_{R}}h(v-u)\text{d}x}.
\end{align}
{If $f$ in decreasing in $t>0$, we infer from $(F_3)$}, $(F_1)$ and $1+\theta_0>0$ 
that
\begin{align}
\int_{B_{R}}
|F(v^+)-F(u^+)|\text{d}x
 =& \int_{B_{R}\cap \{v^+\geq u^+\}}
(F(v^+)-F(u^+))\text{d}x+\int_{B_{R}\cap \{v^+<u^+\}}
(F(u^+)-F(v^+))\text{d}x\notag\\
\leq & \int_{B_{R}\cap \{v^+\geq u^+\}}F(v^+-u^+)\text{d}x+\int_{B_{R}\cap \{v^+<u^+\}}
F(u^+-v^+)\text{d}x\notag\\
\leq & F(1) \max\bigg\{\int_{B_{R}}
|v-u|^{1+\theta_0}\text{d}x, \int_{B_{R}}|v-u|^{1+f_0}\text{d}x\bigg\}\notag\\
=& F(1)\int_{B_{R}}
|v-u|^{1+\theta_0}\text{d}x\notag\\
\leq & C(1+G(1))\int_{B_{R}}
|v-u|^{\gamma}\text{d}x\label{31}
\end{align}
where {}{$0<\gamma<\min\{1,1+\theta_0\} $}, and without loss of generality we assume that $\|v-u\|_{L^{\infty}(B_R)}\leq 1$ due to the boundedness of $v$ and $u$.

By \cite[(10)-(15) on page 44]{ZFZ}, we have
\begin{align}\label{5273}
\int_{B_{R}}
|v-u|^{\gamma}
 \leq& C(\varepsilon_0)R^{n+\alpha_0}+\varepsilon_0 R^{\beta_0}\int_{B_R} G(|\nabla v-\nabla u|)\text{d}x,
\end{align}
where $\alpha_0,\beta_0>0$ are independent of $R$, $\varepsilon_0$ will be chosen later.

If $f$ in increasing in $t>0$, we infer from $(F_4)$, $(f_1)$, $(f_2)$ and the boundedness of $v$ and $u$ that
\begin{align}\label{5273'}
\int_{B_{R}}
|F(v^+)-F(u^+)|\text{d}x=&\int_{B_{R}\cap \{v^+=u^+=0\}}
|F(v^+)-F(u^+)|\text{d}x+\int_{B_{R}\cap \{v^++u^+\neq0\}}
|F(v^+)-F(u^+)|\text{d}x\notag\\
 \leq& \int_{B_{R}\cap \{v^++u^+\neq0\}}
f(\xi)|v^+-u^+|\text{d}x\notag\\
 \leq& \int_{B_{R}\cap \{v^++u^+\neq0\}}
f(\|\xi\|_{L^\infty(B_R)})|v^+-u^+|\text{d}x\notag\\
 \leq& \int_{B_{R}\cap \{v^++u^+\neq0\}}
f(\|u\|_{L^\infty(\Omega)})|v^+-u^+|\text{d}x\notag\\
 \leq& \int_{B_{R}\cap \{v^++u^+\neq0\}}
\frac{1+f_0}{1+\theta_0}\max\{\|u\|_{L^\infty(\Omega)}^{\theta_0},\|u\|_{L^\infty(\Omega)}^{f_0}\}f(1)|v^+-u^+|\text{d}x\notag\\
&\leq  C (1+g(1))\int_{B_{R}}
|v-u|\text{d}x\notag\\
&\leq  C (1+G(1))\int_{B_{R}}
|v-u|\text{d}x,
\end{align}
where $\xi \in (\min\{u^+,v^+\}, \max\{u^+,v^+\})\subset (0, \|u\|_{L^\infty(\Omega)})$ and $C $ is independent of $R$.

Similarly, we get by \cite[(16) on page 44]{ZFZ}
\begin{align}\label{5273''}
\int_{B_{R}}
|v-u|\text{d}x
 \leq& C(\varepsilon_1)R^{n+\alpha_1}+\varepsilon_1 R^{\beta_1}\int_{B_R} G(|\nabla v-\nabla u|)\text{d}x,
\end{align}
where $\alpha_1,\beta_1>0$ are independent of $R$, $\varepsilon_1$ will be chosen later.

{}{Now we estimate $\int_{B_{R}}h(v-u)\text{d}x$. By H\"{o}lder's inequality, Sobolev embedding theorem, and $(\tilde{G}_2)$, it follows
\begin{align}
\int_{B_{R}}|h(v-u)|\text{d}x \leq & \left(\int_{B_{R}} |h|^N \text{d}x\right)^{\frac{1}{N}} \left( \int_{B_{R}} |v-u|^{\frac{N}{N-1}} \text{d}x\right)^{\frac{N-1}{N}}\notag\\
\leq & |B_R|^{\frac{1}{N}-\frac{1}{m}}\left(\int_{B_{R}} |h|^m \text{d}x\right)^{\frac{1}{m}} \left( \int_{B_{R}} |v-u|^{\frac{N}{N-1}} \text{d}x\right)^{\frac{N-1}{N}}\notag\\
\leq & |B_R|^{\frac{1}{N}-\frac{1}{m}}\left(\int_{B_{R}} |h|^m \text{d}x\right)^{\frac{1}{m}}  \left( \int_{B_{R}} |v-u|^{\frac{N}{N-1}} \text{d}x \right)^{\frac{N-1}{N}}\notag\\
\leq & |B_R|^{\frac{1}{N}-\frac{1}{m}}\|h\|_{L^{m}(B_R)}\int_{\Omega} |\nabla (v-u)| \text{d}x \notag\\
\leq & |B_R|^{\frac{1}{N}-\frac{1}{m}} \left( \varepsilon_2 \int_{B_{R}} G(|\nabla u-\nabla v|) \text{d}x +C(\varepsilon_2)|B_{R}|\right),\notag\\
=&\varepsilon_2|B_R|^{\frac{1}{N}-\frac{1}{m}}   \int_{B_{R}} G(|\nabla u-\nabla v|) \text{d}x +C(\varepsilon_2)|B_{R}|^{1+\frac{1}{N}-\frac{1}{m}},\label{1115}
\end{align}
where $\varepsilon_2>0$ will be chosen later.}
{}{By Lemma \ref{Lemma +2.4}, \eqref{5281}, \eqref{31}, \eqref{5273}, \eqref{5273'}, \eqref{5273''}, and \eqref{1115}, we always have
\begin{align*}
\int_{B_{R}} G(|\nabla u-\nabla v|)\text{d}x\leq & CR^{N+\alpha_2}+C\varepsilon R^{\beta_2}\int_{B_{R}}G(|\nabla u-\nabla v|)\text{d}x+CR^{\frac{\lambda}{2}+\frac{N+\alpha_2}{2}}\notag\\
& +C \varepsilon^{\frac{1}{2}}R^{\frac{\lambda}{2}+ \frac{\beta_2}{2}}\bigg(\int_{B_{R}}G(|\nabla u-\nabla v|)\text{d}x\bigg)^{\frac{1}{2}}\notag\\
\leq &CR^{N+\alpha_2}+C\varepsilon R^{\beta_2}\int_{B_{R}}G(|\nabla u-\nabla v|)\text{d}x+CR^{\frac{\lambda+N+\alpha_2}{2}}+C R^{ \lambda+\beta_2}\notag\\
&+ \varepsilon \int_{B_{R}}G(|\nabla u-\nabla v|)\text{d}x,
\end{align*}
where $\alpha_2=\min\{\alpha_0,\alpha_1,1-\frac{N}{m}\},\beta_2=\min\{\beta_0,\beta_1,1-\frac{N}{m}\}, \varepsilon=\max\{\varepsilon_0,\varepsilon_1,\varepsilon_2\}$.}

Choosing $\varepsilon$ small enough, we get
\begin{align}
\int_{B_{R}} G(|\nabla u-\nabla v|)\text{d}x\leq & CR^{s},\label{4.231}
\end{align}
where $s=\min\{ N+\alpha_2,\frac{\lambda+n+\alpha_2}{2},\lambda+\beta_2\}$.

Finally, we get by \eqref{5281}, \eqref{31}, \eqref{5273}, \eqref{5273'}, \eqref{5273''}, \eqref{1115}, and \eqref{4.231}
\begin{align}
\int_{B_{R}} (G(|\nabla u|)-G(|\nabla v|))\text{d}x\leq
CR^{N+\alpha_2}+CR^{\beta_2+s}.\label{4.233}
\end{align}
Putting \eqref{4.233} into \eqref{4.11}, we obtain for all $0<r\leq R$
\begin{align*}
\int_{B_{r}}\! G(|\nabla u-(\nabla u)_{r}|)\text{d}x\leq &
C\bigg(\frac{r}{R}\bigg)^{N+\sigma}\!\!\!\!\int_{B_{R}}\! G(|\nabla u-(\nabla
u)_{R}|)\text{d}x+CR^{N+\alpha_2}+CR^{\beta_2+s}+ CR^{\frac{\lambda}{2}+\frac{N+\alpha_2}{2}}+ CR^{\frac{\lambda}{2}+\frac{\beta_2+s}{2}}.
\end{align*}
Due to the arbitrariness of $\lambda\in (0,N)$, we get $\min\{\beta_2+s,  \frac{\lambda}{2}+\frac{N+\alpha_2}{2},\frac{\lambda}{2}+\frac{\beta_2+s}{2}\}>N$ by setting $\min\{\ \lambda+\alpha_2,\lambda+\beta_2\}>N$. We conclude that there exists $\alpha_3>0$ such that
\begin{align*}
\int_{B_{r}} G(|\nabla u-(\nabla u)_{r}|)\text{d}x\leq &
C \bigg(\frac{r}{R}\bigg)^{N+\sigma}\int_{B_{R}} G(|\nabla u-(\nabla
u)_{R}|)\text{d}x+CR^{N+\alpha_3}.
\end{align*}
In view of Lemma \ref{Lemma 2.13}, we conclude that
there is a constant $\alpha_4\in (0,1)$ such that
\begin{align}
\int_{B_{r}} G(|\nabla u-(\nabla u)_{r}|)\text{d}x\leq Cr^{N+\alpha_4}.\label{4.24}
\end{align}
Proceeding exactly as in \cite[(22) on page 46)]{ZFZ}, we conclude that there is a constant $\alpha\in
(0,1)$ such that
\begin{align}
\int_{B_{r}} |\nabla u-(\nabla u)_{r}|\text{d}x\leq Cr^{N+\alpha},\label{4.25}
\end{align}
which and Campanato's Embedding Theorem
give the H\"{o}lder continuity of the gradient of
$u$.\hfill $\blacksquare$
 \end{pf*}
 \begin{pf*}{Proof of Theorem \ref{Theorem Log-Lip}}
For any fixed $x_{0}\in \Omega$, let $R>0$ such that
$R<\text{dist}(x_{0},\partial \Omega)$. As before, we denote
$B_{R}=B_{R}(x_{0})$.  Let $h$ be the $G-$harmonic function in $B_{R}$ that
agrees
 with $u$ on the boundary, i.e.,
\begin{align}
\text{div}\ \frac{g(|\nabla h|)}{|\nabla h|}\nabla h=0\ \text{in}\
B_{R}\ \ \text{and}\ \ h-u\in W^{1,G}_{0}(B_{R}).\notag
\end{align}
It suffices to note that
$
\int_{B_{R}} (\lambda_{+}\chi_{\{h>0\}}-\lambda_{+}\chi_{\{u>0\}})\text{d}x\leq
\lambda_{+}R^N,
$ {}{and \eqref{1115} becomes
\begin{align*}
\int_{B_{R}}|h(v-u)|\text{d}x \leq\varepsilon_2 \int_{B_{R}} G(|\nabla u-\nabla v|) \text{d}x +C(\varepsilon_2)|B_{R}|.
\end{align*}}
Proceeding as in the proof of Theorem \ref{Theorem +1.1}, we have
\begin{align*}
\int_{B_{r}} G(|\nabla u-(\nabla u)_{r}|)\text{d}x\leq &
C \bigg(\frac{r}{R}\bigg)^{N+\sigma}\int_{B_{R}} G(|\nabla u-(\nabla
u)_{R}|)\text{d}x+CR^N,
\end{align*}
for all $0<r\leq R$. Then Lemma \ref{Lemma 2.13} gives
\begin{align*}
\int_{B_{r}} G(|\nabla u-(\nabla u)_{r}|)\text{d}x\leq Cr^{N}.\notag
\end{align*}
Finally,
\begin{align*}
\int_{B_{r}} |\nabla u-(\nabla u)_{r}|\text{d}x\leq Cr^{N},\notag
\end{align*}
which shows that the gradient of $u$ lies in BMO space and for any
fixed subdomain $\Omega'\Subset\Omega$, there holds
$
\|u\|_{BMO(\Omega')}\leq
C
$ for a universal constant $C>0$.
The residual argument is the same as in \cite[Section 5]{LdT}, and
the desired result can be obtained.\hfill $\blacksquare$
\end{pf*}

 \section{Growth rates near the free boundary for nonnegative minimizers of $ \mathcal {J}(u)$}\label{Sec: 5}
    In view of Corollary \ref{Corollary 17} or Remark \ref{Remark 4}, we may consider non-negative minimizers of $ \mathcal {J}(u)$ and establish their growth rates near the free boundary $\partial \{u>0\}$ for $\lambda_+=0$ and $\lambda_+>0$ respectively.
    {}{To do this, we always assume that $\theta_0<\delta_0$ and \eqref{1115+1} holds. Moreover, we assume that there exists $\tau\in (0,1]$}
such that
\begin{align}\label{+032}
\int_{t}^{t+k}|Q'(s)|\text{d}s\leqslant
c_{0}\left(\frac{k}{t}\right)^{\tau},
\end{align}
for all $t>0,k>0$, where $Q(s)=\frac{tg'(t)}{g(t)}$, $c_{0}=c_{0}({\delta_0},g_{0},\tau)$ is a positive constant.

\begin{theorem}[Growth rates for $\lambda_+=0$]\label{Theorem 1.2}
 {}{Assume that $ h \in L^{m}(\Omega)$ with $m>N$.} Let $u$ be a non-negative minimizer of $ \mathcal {J}(u)$ with $\lambda_+=0$, and $x_0\in \partial\{u>0\}$, $B_{r_0}(x_0)\Subset \Omega$. Then there exists universal constants
$C_0,C_1$, depending only on $N,\theta_0,f_0,\delta_0,g_{0},G(1),\frac{1}{G(1)}$, $\|h\|_{L^{m}(\Omega)}$, $\|q\|_{L^{\infty}(\Omega)},\|\varphi
\|_{L^{\infty}(\Omega)},$ , $\|\varphi
\|_{W^{1,G}(\Omega)}$ and $B_{r_0}(x_0)$, such that
\begin{align}
&|u(x)|\leq C_0{\Phi}(|x-x_0|),\ \ \ \ \forall\ x\in B_{r_0}(x_0),\label{growth1}\\
&|\nabla u(x)|\leq C_1{\Phi}'(|x-x_0|),\ \ \ \ \forall\ x\in B_{r_0}(x_0).\label{growth2}
\end{align}
 for all $0<r<r_0$, where $\Phi(t)=t^{p_0}$ with $ p_0=\min\{\frac{1+g_0}{g_0-\theta_0},\frac{1+\delta_0}{\delta_0-\theta_0},\frac{1+g_0}{g_0}\}>1$.
 \end{theorem}
\begin{theorem}[Growth rates for $\lambda_+>0$]\label{Theorem 1.2'}
 {}{Assume that $ h \in L^{N}(\Omega)$.} Let $u$ be a non-negative minimizer of $ \mathcal {J}(u)$ with $\lambda_+>0$. Assume that $x_0\in \partial\{u>0\}$ and $B_{r_0}(x_0)\Subset \Omega$. Then there exists a universal constant
$C_2$, depending only on $N,$ $\theta_0$, $f_0$, $\delta_0$, $g_{0}$, $G(1)$, $\frac{1}{G(1)}$, $\|h\|_{L^{m}(\Omega)}$, $\|q\|_{L^{\infty}(\Omega)},\|\varphi
\|_{L^{\infty}(\Omega)}$, $\|\varphi
\|_{W^{1,G}(\Omega)}$ and $B_{r_0}(x_0)$ such that
\begin{align}
|u(x)|\leq C_2|x-x_0|,\ \ \ \ \forall\ x\in B_{r_0}(x_0),\label{growth1'}
\end{align}
 for all $0<r<r_0$.
 \end{theorem}
\begin{pf*}{Proof of Theorem \ref{Theorem 1.2}} Due to the local property, we may assume that $u$ is a non-negative minimizer of $\mathcal
{J}(u)$ associated with the domain $B_1(x_0)$ with $x_0 = 0$. Firstly, we prove \eqref{growth1}. Let $
S(j,u) =\sup\limits_{x\in B_{2^{-j}}}|u(x)|.\notag
$
It suffices to show that for all $j\in \mathbb{N}$ there holds
\begin{align}\label{2.1}
S(j+1,u) \leq
\max\left\{c {\Phi}(2^{-j}),S(j,u){\Phi}(2^{-1}),...,
S(j-m,u){\Phi}(2^{-(m+1)}),...,
S(0,u){\Phi}(2^{-j-1})\right\},
\end{align}
with some constant $c>0$.
We prove by contradiction. Let us
suppose \eqref{2.1} fails. Then for any $k\in \mathbb{N}$, there
exists a
sequence of integers $j_k \in \mathbb{N}$ and a sequence of minimizers $u_k$ such that
\begin{align}\label{2.2}
S(j_k+1,u_k) >
\max\left\{k{\Phi}(2^{-j_k}),S(j_k,u_k){\Phi}(2^{-1}),...,
S(j_k-m,u_k){\Phi}(2^{-(m+1)}),...,
S(0,u_k){\Phi}(2^{-j_k-1})\right\}.
\end{align}
Notice that by \eqref{2.2} and the boundedness of $u_k$, it
follows that $j_k\rightarrow \infty$ as $k\rightarrow \infty$.

 Let $v_{k}(x)=\frac{u_k(2^{-j_k}x)}{S(j_k+1,u_k)}$,  $\sigma_k=2^{j_k}S(j_k+1,u_k)$, $G_{k}(t)=\frac{G(\sigma_kt)}{\sigma_kg(\sigma_k)}$ with $g_k(t)=G_k'(t)$, $F_{k}(t)=\frac{F(S(j_k+1,u_k)t)}{\sigma_jg(\sigma_k)}$ with $f_k(t)=F_k'(t)$, {$q_k(x)=q(2^{-j_k}x)$} and $h_k(x)=\frac{S(j_k+1,u_k) }{\sigma_kg(\sigma_k)}h(2^{-j_k}x)$. By $(G_5)$, $G_k$ and $F_k$ satisfy \eqref{+031} and \eqref{hip-F} with the same constants $\delta_0,g_0,\theta_0$ and$f_0$. For all $k>0$, $v_{k}$ is a
minimizer of the functional
$\int_{B_{2^{j_k}}}\big(G_k(|\nabla
v|)+q_kF_k(v^{+})+h_kv\big)\text{d}x.\notag
$ Indeed, by a simple calculation we have
\begin{align}
\int_{B_{2^{j_k}}}\big(G_k(|\nabla
v_k|)+q_kF_k(v_j^{+})+h_kv_k\big)\text{d}x
=\frac{2^{j_kn}}{\sigma_kg(\sigma_k)}\int_{B_{1}}\big(G(|\nabla
u|)+qF(u^{+})+hu\big)\text{d}x.\notag
\end{align}
Particularly, $v_{k}$ is a
minimizer of the following functional
\begin{align}
\mathcal{J}_{k}=\int_{B_{R}}\big(G_k(|\nabla
v|)+q_kF_k(v^{+})+h_kv\big)\text{d}x,\notag
\end{align}
provided $R=2^{m}<2^{j_k}, m$ is fixed.

Notice that by \eqref{2.2} and the definition of $\Phi$, we have $\sup\limits_{B_R}|v_k|\leq C\Phi(R)$, and
\begin{align*}
S(j_k+1,u_k)\geq k\Phi(2^{-j_k})= k (2^{-j_k})^{p_0},
\end{align*}
which gives
\begin{align}
2^{(1+\delta_0)j_k}(S(j_k+1,u_k))^{\delta_0-\theta_0}\geq & k^{\delta_0-\theta_0}2^{(1+\delta_0)j_k}(2^{-j_k})^{p_0(\delta_0-\theta_0) }\notag\\
=& (2^{j_k})^{1+\delta_0-p_0(\delta_0-\theta_0) }k^{\delta_0-\theta_0}
\notag\\
\geq &k^{\delta_0-\theta_0},\label{42}
\end{align}
and
\begin{align}
2^{(1+\delta_0)j_k}(S(j_k+1,u_k))^{g_0-\theta_0}\geq & k^{g_0-\theta_0}2^{(1+g_0)j_k}(2^{-j_k})^{p_0(g_0-\theta_0) }\notag\\
\geq & k^{g_0-\theta_0},\label{43}
\end{align}
Then we get by $(F_1)$, $(G_2)$, $(G_3)$, \eqref{42}, \eqref{43} and $\sup\limits_{B_R}|v_k|\leq C\Phi(R)$
\begin{align}
|q_kF_k(v_k^{+})|=&|q_k|\frac{F(S(j_k+1,u_k) v_k^{+})}{\sigma_kg(\sigma_k)}\notag\\
\leq &  \frac{C|S(j_k+1,u_k)|^{1+\theta_0}F( v_k^{+})}{\sigma_kg(\sigma_k)}\notag\\
\leq &  \frac{C|S(j_k+1,u_k)|^{1+\theta_0}}{\min\{(2^{j_k}S(j_k+1,u_k))^{1+\delta_0}, (2^{j_k}S(j_k+1,u_j))^{1+g_0}\}}\notag\\
\leq &  \frac{C}{\min\{2^{(1+\delta_0)j_k}(S(j_k+1,u_k))^{\delta_0-\theta_0}, (2^{(1+g_0)j_k}(S(j_k+1,u_k))^{g_0-\theta_0}\}}\notag\\
\leq &  \frac{C}{\min\{ k^{\delta_0-\theta_0},k^{g_0-\theta_0}\}}\rightarrow 0,\ \ \text{as}\ k\rightarrow \infty.
\end{align}
Similarly, due to that $
2^{(1+\delta_0)j_k}(S(j_k+1,u_k))^{\delta_0}\geq k^{\delta_0}$, and $2^{(1+\delta_0)j_k}(S(j_k+1,u_k))^{g_0}\geq k^{g_0},
$ we have {}{$\|h_k\|_{L^{m}(B_R)}\rightarrow 0\ \ \text{as}\ k\rightarrow \infty$}.
For $k$ large enough,
according to the $C^{1,\alpha}$ regularity of minimizers, we obtain $
\|v_{k}\|_{C^{1,\alpha}(B_{R})}\leq C
$ (see Theorem \ref{Theorem +1.1}). {Note that $C$ depends on $\frac{1}{G_k(1)}$ and $ G_k(1)$. However by $(G_5)$, we see that $C$ depends on $\frac{1}{G(1)}$ and $ G(1)$ essentially, thus it is independent of $k$.} Therefore, up to subsequence, we get $v_{k}\rightarrow v_{0}$
in
$C^{1,\beta}(B_{r_0})$ with $0<\beta<\alpha$ and any $r_0< 1$. We deduce by $v_{k}(0)=0$ and $
\sup\limits_{B_{\frac{1}{2}}}|v_{k}|=1 $ that
\begin{align} \label{contradiction}\sup_{B_{\frac{1}{2}}}|v_{0}|=1,v_{0}(0)=0,
\end{align}
On the other hand, using the compact condition \eqref{+032}, we conclude that (see \cite[Theorem 6.1]{BM}) there exists a function $G_{\infty}\in C^{2}(0,+\infty)$ such that, up to a subsequence,
\begin{align*}
&G_{k}\rightarrow G_{\infty},\ \ g_k=G'_{k}\rightarrow G'_{\infty}=g_\infty\ \  \text{uniformly\  in compact \ subsets\  of}\ [0,+\infty),\\
&G''_{k}\rightarrow G''_{\infty}\ \ \text{uniformly\  in compact \ subsets\  of}\ (0,+\infty),
\end{align*}
and $g_{\infty}$ satisfies structural condition \eqref{+031} with the same constants. Furthermore, we infer that $v_0$ is a $G_\infty-$harmonic function in $B_{1}$.
 Since $v_k\geq 0$ in $B_{1}$, $v_0\geq 0$ in $B_{r_0}$. Recalling $v_0(0)=0$ and the Harnack's inequality, we have $v_0\equiv 0$ in $B_{r_0}$. Finally we get $v_0\equiv 0$ in $B_{1}$ due to the continuity of $v_0$ and the arbitrariness of $r_0$.
which is a contradiction with \eqref{contradiction}. Therefore have proved \eqref{growth1}.

 Now, we prove \eqref{growth2}. Set
$
S(j,|\nabla u|) =\sup\limits_{x\in B_{2^{-j}}}|\nabla u(x)|.\notag
$
It suffices to show
\begin{align}\label{4.6'}
S(j+1,|\nabla u|) \leq
\max\bigg\{c {\Phi}'(2^{-j}),S(j,|\nabla u|){\Phi}'(2^{-1}),&...,
S(j-m,|\nabla u|){\Phi}'(2^{-(m+1)}),\notag\\
&...,
S(0,|\nabla u|){\Phi}'(2^{-j-1})\bigg\}.
\end{align}
for some positive constant $c$. By contradiction, suppose that
\eqref{4.6'} fails. Then for any $k\in \mathbb{N}$, there exists a
sequence of integers $j_k$ and a sequence of minimizers
$u_k$ such that
\begin{align}\label{1.6-}
S(j_k+1,|\nabla u_k|) >
\max\big\{k {\Phi}'(2^{-j_k}),S(j_k,|\nabla u_k|){\Phi}'(2^{-1}),&...,
S(j_k-m,|\nabla u_k|){\Phi}'(2^{-(m+1)}),\notag\\
&...,
S(0,|\nabla u_k|){\Phi}'(2^{-j_k-1})\big\}.
\end{align}
Let
$v_{k}(x)=\frac{u_k(2^{-j_k}x)}{2^{-j_k}S(j_k+1,u_k)}$,  $\varrho_k=S(j_k+1,u_k)$, $G_{k}(t)=\frac{G(\varrho_kt)}{\varrho_kg(\varrho_k)}$ with $g_k(t)=G_k'(t)$, $F_{k}(t)=\frac{F(S(j_k+1,u_k)t)}{\sigma_jg(\varrho_k)}$ with $f_k(t)=F_k'(t)$, {$q_k(x)=q(2^{-j_k}x)$} and $h_k(x)=\frac{S(j_k+1,u_k) }{\varrho_kg(\varrho_k)}h(2^{-j_k}x)$. Then for all $k>0$, $v_{k}$ is a
minimizer of the functional
$\int_{B_{2^{j_k}}}\big(G_k(|\nabla
v|)+q_kF_k(v^{+})+h_kv\big)\text{d}x.\notag
$ By \eqref{growth1} and \eqref{1.6-}, we have $\sup\limits_{B_{1}}|v_{k}|\leq \frac{C}{k}\rightarrow 0$ as $k \rightarrow \infty$. Arguing as before, we get
$
|q_kF_k(v_k^{+})|\rightarrow 0\ \ \text{and}\ \ {}{\|h_k\|_{L^m(B_R)}}\rightarrow 0 \ \text{as}\ k\rightarrow \infty
$, and we can conclude that there exists a $G_\infty-$harmonic function $v_0$ in $B_{1}$, satisfying $v_{k}\rightarrow v_{0}$
in
$C^{1,\beta}(B_{r_0})$ with some $\beta\in (0,1)$ and any $r_0< 1$. Furthermore, we conclude that $v_{0}\equiv0$
in $B_{1}$. However, note that $ \sup\limits_{B_{1/2}}|\nabla
v_{j_{k}}|=1.\notag $ Thus $ \sup\limits_{B_{1/2}}|\nabla v_{0}|=1,\notag $
which gives a contradiction.\hfill $\blacksquare$
\end{pf*}
\begin{pf*}{Proof of Theorem \ref{Theorem 1.2'}}
Let $\Phi(t)=t^{p_0}$ for all $t\geq 0$, where $p_0=\min\{\frac{1+g_0}{g_0-\theta_0},\frac{1+\delta_0}{\delta_0-\theta_0},\frac{1+g_0}{g_0},1\}=1$. Then one can proceed with a slight modification of the proof of Theorem \ref{Theorem 1.2} to obtain $ |u(x)|\leq C_2 \Phi(|x-x_0|)$.
\hfill $\blacksquare$
\end{pf*}
\begin{corollary}[Optimal growths in the non-homogenous one-phase problems for $p-$Laplacian]
  Let $G(t)=t^p$ with $p>1$, and $F(t)=t^{\gamma}$ with $0<\gamma<p$. Let $u$ be a nonnegative minimizer of $ \mathcal {J}(u)$ with $\lambda+=0$ in \eqref{main equation} and $x_0\in \partial\{u>0\}$. Then there exists a universal constant
$C$ such that
\begin{align*}
|u(x)|\leq C|x-x_0|^{p_0},\ |\nabla u(x)|\leq C|x-x_0|^{p_0-1},\ \ \ \forall\ x\in B_{r_0}(x_0)\Subset \Omega
\end{align*}
 for all $0<r<r_0$, where $p_0=\min\{\frac{p}{p-\gamma},\frac{p}{p-1}\}>1$.
 \end{corollary}
\begin{remark} Checking the proof of Theorem \ref{Theorem 1.2}, if $h=0$ and $u$ is a nonnegative minimizer of $ \mathcal {J}(u)$ with $\lambda+=0$ in \eqref{main equation}, then we have
 \begin{align*}
|u(x)|\leq C|x-x_0|^{p_1},\ |\nabla u(x)|\leq C|x-x_0|^{p_1-1},\ \ \ \forall\ x\in B_{r_0}(x_0)\Subset \Omega
\end{align*}
where $p_1=\min\{\frac{1+g_0}{g_0-\theta_0},\frac{1+\delta_0}{\delta_0-\theta_0}\}$. Particularly, if $G(t)=t^p,p>1$ and $F(t)=t^{\gamma},0<\gamma<p$, we have
\begin{align*}
|u(x)|\leq C|x-x_0|^{\frac{p}{p-\gamma}},\ |\nabla u(x)|\leq C|x-x_0|^{\frac{\gamma}{p-\gamma}},\ \ \ \forall\ x\in B_{r_0}(x_0)\Subset \Omega,
\end{align*}
which are the optimal growth rates of minimizers and their gradients in the homogeneous one-phase free boundary problems for $p-$Laplacian.
\end{remark}
\begin{remark}
 Condition \eqref{+032} is used only for the compactness of $G_k$ by blow-up techniques, see, e.g., \cite[Theorem 6.1]{BM}. For the case of $p-$Laplacian, i.e., $G(t)=t^p$, \eqref{+032} becomes trivial due to that $ Q'(s)\equiv 0$.
 \end{remark}
\section{Local Lipschitz continuity of non-negative minimizers of $ \mathcal {J}(u)$ with $\lambda_+>0$}\label{Sec: Lipschitz continuity}
In this section, in order to obtain local Lipschitz continuity of non-negative minimizers of $ \mathcal {J}(u)$ with $\lambda_+>0$, we make further assumptions on $F$, i.e., assume that $F\in C^1([0,+\infty); [0,+\infty))$. Note that $f\in C([0,+\infty); [0,+\infty))$ and there exists positive constants $C_1$ and $C_2$ such that $f(t)\leq C_1+C_2g(t)$ for all $t\geq 0$. {}{We also assume that \eqref{1115+1} and \eqref{+032} hold, and $ h\in L^{\infty}(\Omega)$.}

 We say that a function $u\in W^{1,G}(D)$ is a weak solution of the equation $
\text{div} \ \frac{g(|\nabla u|)}{|\nabla u|}\nabla u=qf(u)+h\ \ \ \text{in}\ D\subset \Omega,
$
if
\begin{align}\label{weak solution}
\int_{D}\frac{g(|\nabla u|)}{|\nabla u|} \nabla u \nabla \xi \text{d}x  +\int_{D} (qf(u)+h) \xi \text{d}x=0
\end{align}
holds for all $\xi\in W_0^{1,\tilde{G}}(D)$, where $D $ is a domain.
\begin{lemma}[Harnack's inequality ]\label{Harnack's inequality}
 Let $u\in W^{1,G}(B_{R})$ with $0\leq u\leq M$
is {a weak solution} of $\text{div} \ \frac{g(|\nabla u|)}{|\nabla u|}\nabla u=qf(u)+h$ in $B_{R}$. Then, there exists a universal constant $t_0>0$ and a constant $C_{r}>0$ depending only on ${\delta_0},g_{0},M,t_0,\|q\|_{L^\infty(B_{R})}$ and $R-r$ such that
\begin{align}
\sup_{B_{r}}u\leq
C_{r}\bigg(\inf_{B_{r}}u+g^{-1}(R)R\bigg),\notag
\end{align}
for all $0<r\leq R$.
 \end{lemma}
 \begin{pf*}{Proof} It is a direct result of \cite[ Corollary 1.4]{L}. Indeed, we amy set $a_1=a_2=a_4=a_5=0$ and $a_3=1$ in (1.3a) and (1.3b) of \cite[ Corollary 1.4]{L} for our problem. We shall verify that conditions (1.3c)'' and (1.4) of \cite[ Corollary 1.4]{L} are satisfied. By $(f_2)$, there exits $t_0>0$ such that for all $t>t_0$, there holds
 \begin{align*}
f(t)\leq &  g(t)=g\bigg(\frac{t}{R}\cdot R\bigg)\notag\\
\leq & g\bigg(\frac{t}{R}\bigg)\frac{t}{R}\frac{R}{t}\max\{R^{{\delta_0}},R^{g_0}\}\ \ \ \ \text{by}\ (g_1)\notag\\
\leq & \frac{1}{t_0}\cdot g\bigg(\frac{t}{R}\bigg)\frac{t}{R},
\end{align*}
 where without loss of generality we assume that $R\leq 1$. For $0\leq t\leq t_0$, due to $F\in C^1([0,+\infty); [0,+\infty))$, $f(t)=F'(t) $ is continuous in $[0,t_0]$. Then there exists a constant $M_0>0$ such that $f(t)\leq M_0$ for all $t\in [0,t_0]$. Then for all $t\geq 0$, we have
 \begin{align*}
f(t)\leq \frac{1}{t_0}\cdot g\bigg(\frac{t}{R}\bigg)\frac{t}{R}+M_0.
\end{align*}
Thus we can choose $b_0=0,b_1= \frac{1}{t_0}\cdot\|q\|_{L^\infty(B_{R})}, b_2=M_0\cdot\|q\|_{L^\infty(B_{R})}+\|h\|_{L^\infty(B_{R})}$ and $\chi=g^{-1}(b_2R)$ in (1.3c)'' and (1.4) of \cite[ Corollary 1.4]{L}. Finally, \text{by}\ $(\tilde{g}_1)$ and \cite[ Corollary 1.4]{L}, we have
\begin{align*}
\sup_{B_{r}}u\leq
C\bigg(\inf_{B_{r}}u+g^{-1}(b_2R)R\bigg)
\leq  C'\bigg(\inf_{B_{r}}u+g^{-1}(R)R\bigg).
\end{align*}\hfill $\blacksquare$
 \end{pf*}
 A consequence of Theorem \ref{Holder continuity} is the fact that $\{u>0\}$ is an open set. We have the following result.
 \begin{lemma}
{Let $u$ be a non-negative minimizer of $ \mathcal {J}(u)$ with $\lambda_+>0$ in \eqref{main equation}. Then $u$ is a weak solution of the following equation
  \begin{align*}
 \text{div}\ \frac{g(|\nabla u|)}{|\nabla u|}\nabla u=qf(u)+h\ \ \ \ \text{in }\ \{u>0\}.
 \end{align*}}
  \end{lemma}
\begin{pf*}{Proof}
For any ball $B\subset\{u>0\} $, consider first that $\xi \in C^{\infty}_{0}(B)$. There exists $0<\varepsilon_0 \leq 1$ small enough such that $\{u \pm\varepsilon \xi> 0\}\cap B =B$ for all $0<\varepsilon \leq \varepsilon_0 .$
Standard arguments implies that
\begin{align}\label{lebesgue-nonlinear}
\displaystyle\lim_{\varepsilon \rightarrow 0^{+}} \int_{B}\frac{F((u+ \varepsilon \xi)^+) - F(u^{+})}{\varepsilon} \text{d}x=  \int_{B} f(u) \xi\text{d}x.
\end{align}
The minimality  of $u$ implies that
\begin{align}\label{subsol}
\begin{aligned}
0 \leq & \frac{1}{\varepsilon}\int_{B} \bigg(G(|\nabla (u+ \varepsilon \xi)|) -  G(|\nabla u|) + q(F((u+\varepsilon \xi)^{+})  - F(u^{+}))\bigg)\text{d} x\\
\leq &\int_{B} g(|\nabla u +\varepsilon \nabla\xi|)\frac{\nabla u +\varepsilon \nabla\xi}{|\nabla u +\varepsilon \nabla\xi|}\nabla \xi\text{d}x + \frac{1}{\varepsilon}\left( \int_{B}q(F((u+\varepsilon \xi)^{+})  - F(u^{+}))\text{d}x+ \int_{B}h\varepsilon \xi\text{d}x\right),
\end{aligned}
\end{align}
where in the last inequality we used the convexity of $G.$

From \eqref{lebesgue-nonlinear}, \eqref{subsol} and letting $\varepsilon \rightarrow 0^{+}$, we get
\begin{align}\label{subsolution}
\int_{B}\frac{g(|\nabla u|)}{|\nabla u|} \nabla u \nabla \xi \text{d}x  +\int_{B} (qf(u)+h) \xi \text{d}x\geq 0.
\end{align}
Using the function  $\phi = u-\varepsilon \xi$ and repeating the previous arguments we get
\begin{align}\label{subsolution2}
-\int_{B}\frac{g(|\nabla u|)}{|\nabla u|} \nabla u \nabla \xi \text{d}x  -\int_{B}  (qf(u)+h) \xi \text{d}x\geq 0,
\end{align}
for all $\xi \in C^{\infty}_{0}(B).$ By \eqref{subsolution} and \eqref{subsolution2}, \eqref{weak solution} holds for all $\xi \in C^{\infty}_{0}(B).$ Now for $\xi \in W_0^{1,\tilde{G}}(B),$ let $\xi_n \in C_0^{\infty}(B)$ with $ \xi_n\rightarrow \xi$ in $ W_0^{1,\tilde{G}}(B)$ as $n\rightarrow \infty$, then \eqref{weak solution} holds with $\xi_n \in C_0^{\infty}(B)$. We conclude the desired result by letting $n\rightarrow \infty$ and the arbitrariness of $B$.\hfill $\blacksquare$
\end{pf*}
\begin{theorem}[Local Lipschitz continuity for $\lambda_+>0$]\label{Lipschitz continuity}
 Given a subdomain $\Omega'\Subset \Omega$, there exists a constant $C > 0$ that
depends only on
$\Omega'$ and universal constants, such that for any nonnegative minimizer of $ \mathcal {J}(u)$ with $\lambda_+>0$ in \eqref{main equation}, there holds
\begin{align}\label{Lipschitz}
\|\nabla u\|_{L^{\infty}(\Omega')}\leq C.
\end{align}
 \end{theorem}
 \begin{pf*}{Proof}
 We proceed as the proof of \cite[Theorem 4.1]{LT}, supposing that \eqref{Lipschitz} fails. Then there exists a sequence of points $X_j \in \Omega'\cap \{u>0\}$ such that
\begin{align} \label{036}
 X_j \rightarrow \partial \{u>0\}\ \ \ \ \text{and}\ \ \ \ \frac{u(X_j)}{\text{dist}(X_j,\partial\{u>0\})}\rightarrow +\infty\ \ \ \ \text{as}\ \ j\rightarrow +\infty.
 \end{align}
 Denote $U_j=u(X_j)$ and $d_j=\text{dist}(X_j,\partial\{u>0\})$. Let $Y_j\in \partial\{u>0\}$ satisfying $d_j=|X_j-Y_j|$.
{Note that we have}
 \begin{align*}
 \text{div}\ \frac{g(|\nabla u|)}{|\nabla u|}\nabla u=qf(u)+h\ \ \ \ \text{in }\ \{u>0\}.
 \end{align*}
Thus, by Harnack's inequality, $(\tilde{g}_{1})$, and the boundedness of $u$, there exists a constant $c$ depending only on $\Omega'$ and universal constants, such that
\begin{align}
d_j+
\inf_{\overline{B}_{\frac{3}{4}d_j}(X_j)}u\geq cU_j.\notag
\end{align}
In turn, we have
\begin{align}\label{147}
\sup_{\overline{B}_{\frac{d_j}{4}}(Y_j)}u\geq cU_j-d_j.
\end{align}
Consider the set
$
A_j=\big\{Z\in B_{d_j}(Y_j): \text{dist}(Z,\partial \{u>0\})\leq \frac{1}{3}\text{dist}(Z,\partial B_{d_j}(Y_j))\big\}.
$
Then $ B_{\frac{d_j}{4}}(Y_j)\subset A_j$ (see the proof of \cite[ Theorem 4.1]{LT}).
Thus
\begin{align*}
\text{dist}(Z_j,\partial B_{d_j}(Y_j))u(Z_j)&:=M_j\notag\\
&:=\sup_{Z\in A_j} \text{dist}(Z,\partial B_{d_j}(Y_j))u(Z)\notag\\
&\geq \sup_{Z\in\overline{B}_{\frac{d_j}{4}}(Y_j)} \text{dist}(Z,\partial B_{d_j}(Y_j))u(Z)\notag\\
&\geq \sup_{Z\in\overline{B}_{\frac{d_j}{4}}(Y_j)}\frac{3d_j}{4} u(Z)\notag\\
=&\frac{3d_j}{4}\sup_{\overline{B}_{\frac{d_j}{4}}(Y_j)} u.
\end{align*}
It follows that
\begin{align*}
u(Z_j)\geq \frac{d_j}{\text{dist}(Z_j,\partial B_{d_j}(Y_j))}\frac{3}{4} \sup_{\overline{B}_{\frac{d_j}{4}}(Y_j)} u\geq \frac{3}{4}\sup_{\overline{B}_{\frac{d_j}{4}}(Y_j)} u.
\end{align*}
Using \eqref{147}, we have
\begin{align}\label{1444}
u(Z_j)\geq \frac{3}{4}(cU_j-d_j).
\end{align}
For each $j$, let $W_j\in \partial\{u>0\}$ satisfy
\begin{align*}
r_j=|Z_j-W_j|=\text{dist}(Z_j,\partial \{u>0\})\leq \frac{1}{3}\text{dist}(Z_j,\partial B_{d_j}(Y_j)).
\end{align*}
One may get (see (4.7) in \cite{LT})
\begin{align}\label{1447}
\frac{d_j}{r_j}\geq 4.
\end{align}
It follows from \eqref{1444}, \eqref{1447} and \eqref{036} that
\begin{align}\label{1449}
\frac{u(Z_j)}{r_j}\geq \frac{3d_j}{4r_j}\bigg(c\frac{U_j}{d_j}-1\bigg)\geq 3 \bigg (c\frac{U_j}{d_j}-1\bigg)\rightarrow +\infty.
\end{align}
Proceeding as (4.10), (4.11) in \cite{LT}, one has (for $j$ large enough)
\begin{align}\label{14410411}
\sup_{B_{\frac{r_j}{2}}(W_j)}u\leq 2u(Z_j),\ \ \ \ \sup_{\overline{B}_{\frac{r_j}{4}}(W_j)}\frac{u}{u(Z_j)}\geq \frac{c'}{2},
\end{align}
for some universal constant $c'>0$.
Now for each $j$, define the function $u_j:B_{1}(0)\rightarrow (0,2)$ by
\begin{align}\label{14412}
u_j(X)=\frac{u(W_j+\frac{r_j}{2}X)}{u(Z_j)}.
\end{align}
It follows from \eqref{14410411} that
\begin{align}\label{14413}
\max_{B_1(0)}u_j\leq 2,\ \ \ \ \max_{B_1(0)}u_j\geq \frac{c'}{2},\ \ \ \ u_j(0)=0.
\end{align}

 Let $\sigma_j=\frac{2u(Z_j)}{r_j}$, $G_{j}(t)=\frac{G(\sigma_jt)}{\sigma_jg(\sigma_j)}$ with $g_j(t)=G_j'(t)$, $F_{j}(t)=\frac{F(u(Z_j)t)}{\sigma_jg(\sigma_j)}$ with $f_j(t)=F_j'(t)$, {$q_j(X)=q(W_j+\frac{r_j}{2}X)$}, $h_j(X)=\frac{h(W_j+\frac{r_j}{2}X) }{\sigma_jg(\sigma_j)}$ and $\lambda_{+j}=\frac{\lambda_{+}}{\sigma_jg(\sigma_j)}$. Then for all $j>0$,
 \begin{align}\label{1614}
 {\delta_0}\leq \frac{tg_j'(t)}{g_j(t)}\leq g_0, \ \ \ \ \ {1+\theta_0}\leq \frac{tF_j'(t)}{F_j(t)}\leq 1+f_0.
 \end{align}
  Let $v$ be the $G-$harmonic function in $B_{\frac{r_j}{2}}(W_j)$ with the boundary data $u$, i.e.,
\begin{align}
\left\{
  \begin{array}{ll}\nonumber
   \text{div}\ \frac{g(|\nabla v|)}{|\nabla v|}\nabla v &=0\ \ \ \ \text{in }\ B_{\frac{r_j}{2}}(W_j), \\
  \ \ \ \ \ \ \ \ \  \ \ \ \ \ \ \ \ \  \ \ v&=u\ \ \ \ \text{on}\ \partial B_{\frac{r_j}{2}}(W_j).
  \end{array}\notag
\right.
\end{align}
Let $v_j: B_1(0)\rightarrow (0,2)$ be defined by $v_j(X)=\frac{v(W_j+\frac{r_j}{2}X)}{u(Z_j)}$. Then $v_j$ satisfies
\begin{align}
\left\{
  \begin{array}{ll}\nonumber
   \text{div}\ \frac{g_j(|\nabla v_j|)}{|\nabla v_j|}\nabla v_j=&0\ \ \ \ \ \text{in }\ B_{1}(0), \\
    \ \ \ \ \ \ \  \ \ \ \ \ \ \ \ \ \ \ \ \ \ \ \ v_j=&u_j\ \ \ \text{on }\ \partial B_{1}(0).
  \end{array}\notag
\right.
\end{align}
 Let $Y=W_j+\frac{r_j}{2}X$, then
  \begin{align*}
\int_{B_{\frac{r_j}{2}}(W_j)}G(|\nabla u(Y)|)\text{d}Y=\bigg(\frac{r_j}{2}\bigg)^n\int_{B_{1}(0)}G(\sigma_j|\nabla u_j(X)|)\text{d}X.
\end{align*}
It follows
\begin{align*}
\int_{B_{\frac{r_j}{2}}(W_j)}\frac{G(|\nabla u(Y)|)}{\sigma_jg(\sigma_j)}\text{d}Y=&\bigg(\frac{r_j}{2}\bigg)^n
\int_{B_{1}(0)}\frac{G(\sigma_j|\nabla u_j(X)|)}{\sigma_jg(\sigma_j)}\text{d}X\notag\\
=&\bigg(\frac{r_j}{2}\bigg)^n\int_{B_{1}(0)}G_j(|\nabla u_j(X)|)\text{d}X,
\end{align*}
 which gives
 \begin{align}\label{+14}
\int_{B_{1}(0)}G_j(|\nabla u_j|)\text{d}x=\bigg(\frac{r_j}{2}\bigg)^{-n}
\int_{B_{\frac{r_j}{2}}(W_j)}\frac{G(|\nabla u|)}{\sigma_jg(\sigma_j)}\text{d}x.
 \end{align}

 By the minimality of $u$, we have
\begin{align}\label{+15}
\int_{B_{\frac{r_j}{2}}(W_j)}G(|\nabla u|)\text{d}x-\int_{B_{\frac{r_j}{2}}(W_j)}G(|\nabla v|)\text{d}x\leq&
\int_{B_{\frac{r_j}{2}}(W_j)}\bigg(q(F(v^+)-F(u^+))+ h(v-u)\bigg)\text{d}x
\notag\\
&+\lambda_+\int_{B_{\frac{r_j}{2}}(W_j)}(\chi_{\{h>0\}}
-\chi_{\{u>0\}})\text{d}x\notag\\
\leq &  \|q\|_{L^\infty(B_{\frac{r_j}{2}})}\int_{B_{\frac{r_j}{2}}(W_j)}(|F(h)|+|F(u)|)\text{d}x+Cr_j^n
\notag\\
\leq &  Cr_j^n,
\end{align}
where we used the boundedness of $h$ and $u,v$, and the increasing property of $F$ in the last inequality.

 We infer from \eqref{+14} and \eqref{+15}
\begin{align}\label{+16}
\int_{B_{1}(0)}G_j(|\nabla u_j|)\text{d}x-\int_{B_{1}(0)}G_j(|\nabla v_j|)\text{d}x\leq \frac{C}{\sigma_jg(\sigma_j)}\rightarrow 0\ \ \text{by}\ \sigma_j\rightarrow +\infty.
\end{align}
Then we deduce by Lemma \ref{Lemma +2.4} and \eqref{+16}
\begin{align}\label{+17}
\int_{B_{1}}G_j(|\nabla u_j-\nabla v_j|)\text{d}x \leq C\bigg(\frac{1}{\sigma_jg(\sigma_j)}+\frac{1}{\sqrt{\sigma_jg(\sigma_j)}}\bigg)\rightarrow 0\ \ \text{as}\ j\rightarrow +\infty,
\end{align}
where we used the uniform boundedness of $\int_{B_1}G_j(|\nabla v_j|)\text{d}x$ due to the uniform boundedness of $u_j$ and $v_j$ (see, e.g., Lemma \ref{Lemma 2.131}).

 We get by $(G_3)$
 \begin{align}\label{+18}
\int_{B_{1}^-}|\nabla u_j-\nabla v_j|^{1+g_0}\text{d}x+ \int_{B_{1}^+}|\nabla u_j-\nabla v_j|^{1+{\delta_0}}\text{d}x\leq C \int_{B_{1}}G_j(|\nabla u_j-\nabla v_j|)\text{d}x,
\end{align}
 where $B_1^{-}=B_1\cap\{ |\nabla u_j-\nabla v_j|<1\}$ and $B_1^+=B_1\cap\{ |\nabla u_j-\nabla v_j|\geq 1\}$.
H\"{o}lder's inequality gives
\begin{align}\label{+19}
\int_{B_{1}^-}|\nabla u_j-\nabla v_j|^{1+{\delta_0}}\text{d}x\leq C\bigg(\int_{B_{1}^-}|\nabla u_j-\nabla v_j|^{1+g_0}\text{d}x\bigg)^{\frac{1+{\delta_0}}{1+g_0}}.
\end{align}
So we obtain by \eqref{+18} and  \eqref{+19}
\begin{align*}
\bigg(\int_{B_{1}^-}|\nabla u_j-\nabla v_j|^{1+{\delta_0}}\text{d}x\bigg)^{\frac{1+g_0}{1+{\delta_0}}}+ \int_{B_{1}^+}|\nabla u_j-\nabla v_j|^{1+{\delta_0}}\text{d}x\leq C \int_{B_{1}}G_j(|\nabla u_j-\nabla v_j|)\text{d}x,
\end{align*}
which and \eqref{+17} imply that
\begin{align*}
\int_{B_{1}}|\nabla u_j-\nabla v_j|^{1+{\delta_0}}\text{d}x\rightarrow 0\ \ \text{by}\ j\rightarrow +\infty.
\end{align*}
It follows by Poincar\'{e}'s inequality that
\begin{align}\label{+20}
u_j- v_j\rightarrow 0\ \ \text{strongly\ in}\ W^{1,1+{\delta_0}}_{0}(B_1).
\end{align}
Note that the uniform boundedness of $v_j$ guarantees that, for any $r_0\in (0,1)$ there exists a universal constant $C>0$ {satisfying $\|v_j\|_{C^{1,\alpha}(B_{r_0})}\leq C$ (see, e.g. \cite[Theorem 1.2]{ZFZ})}. Therefore, we can find $v_{0}$ in $ B_{r_0}$ such that, up to a subsequence,
\begin{align*}
v_{k}\rightarrow v_{0}\ \text{and}\  \nabla v_{k}\rightarrow \nabla v_{0},\ \text{uniformly\  in}\  B_{r_0}.
\end{align*}
On the other hand, noting that $u_j$ is a minimizer of the following functional
\begin{align}
\mathcal {J}_j=\int_{B_1}(G_j(|\nabla
w|)+q_jF_j(w^+)+h_jw+\lambda_{+j}\chi_{\{w>0\}} )\text{d}x\rightarrow\text{min},\notag
\end{align}
and recalling the structural conditions of $g_j(t)$, $F_j(t)$, and $\lambda_{+j}$, and the boundedness of $q_j,h_j$, we have the uniform H\"{o}lder's estimate of $u_j$, i.e., $\|u_j\|_{C^{\beta}(B_{r_0})}\leq C $. So, we conclude that there exists a $u_{0}\in C^{\beta}(B_{r_0})$ such that
\begin{align*}
u_k\rightarrow u_{0}\ \ \ \text{uniformly\  in}\  B_{r_0}.
\end{align*}
 We conclude this way that $u_{0}=v_{0}$ in $B_{r_0}$ by \eqref{+20}.

 Now using the compact condition \eqref{+031}, we conclude that (see \cite[Theorem 6.1]{BM}) there exists a function $G_{\infty}\in C^{2}(0,+\infty)$ such that, up to a subsequence,
\begin{align*}
&G_{j}\rightarrow G_{\infty},\ \ g_j=G'_{j}\rightarrow G'_{\infty}=g_\infty\ \  \text{uniformly\  in compact \ subsets\  of}\ [0,+\infty),\\
&G''_{j}\rightarrow G''_{\infty}\ \ \text{uniformly\  in compact \ subsets\  of}\ (0,+\infty),
\end{align*}
and $g_{\infty}$ satisfies the same structural condition as \eqref{1614} with the same constants.

 We now claim that $u_0$ is a $G_\infty-$harmonic function in $B_{r_0}$.
  Indeed, by the minimality of $u_j$ again, we have for any $\varphi\in C^{\infty}_0(B_{r_0})$
\begin{align}\label{+21}
\int_{B_{r_0}}G_j(|\nabla
u_j|)\text{d}x
\leq & \int_{B_{r_0}}\bigg(G_j(|\nabla
u_j+\nabla\varphi|)+q_jF_j((u_j+\varphi)^+)-q_jF_j(u_j^+)+h_j((u_j+\varphi)^+-u_j^+)\bigg)\text{d}x\notag\\
&+\int_{B_{r_0}}(\lambda_{+j}\chi_{\{u_j+\varphi>0\}} -\lambda_{+j}\chi_{\{u_j>0\}})\text{d}x.
\end{align}
Note that $q_j$ is uniformly bounded, $\sigma_j\rightarrow +\infty$, $\|h_j \|_{ B_1} \rightarrow 0$ and $u, \varphi$ are bounded, then
\begin{align}
&\|h_j((u_j+\varphi)^+-u_j^+)\|_{ B_1}
\rightarrow0,\label{+22'}\\
&q_jF_j(u_j^+)=\frac{q_jF\big(|u(W_j+\frac{r_j}{2}X)|\big)}{\sigma_jg(\sigma_j)}\leq \frac{\|q\|_{L^{\infty}(B_1)}F\big(\sup\limits_{B_1}|u|\big)}{\sigma_jg(\sigma_j)}
\rightarrow0,\label{+22}
\end{align}
and
\begin{align}\label{+23}
 q_jF_j((u_j+\varphi)^+)=&
\frac{\|q\|_{L^{\infty}(B_1)}F\big(u(W_j+\frac{r_j}{2}X)+u(Z_j)\varphi\big)}{\sigma_jg(\sigma_j)}
\notag\\
\leq & \frac{\|q\|_{L^{\infty}(B_1)}F\big((1+\sup\limits_{B_1}|\varphi|)\sup\limits_{B_1}|u|\big)}{\sigma_jg(\sigma_j)}
\rightarrow0.
\end{align}
Note also that
  \begin{align*}
 G_j(|\nabla u_j|)\leq C(G_j(|\nabla u_j-\nabla v_j|)+G_j(|\nabla v_j|)),
\end{align*}
which, the $C^1-$convergence of $v_j$, and \eqref{+17} imply that there exists $\xi\in L^1(B_{r_0})$ such that
    \begin{align*}
  G_j(|\nabla u_j|)\leq \xi\ \ \ \ \text{a.e.\ in}\ B_{B_{r_0}}.
\end{align*}
Once $\nabla u_j\rightarrow \nabla u_0$ a.e. in $B_{r_0}$, Lebesgue's dominated convergence theorem implies
  \begin{align*}
\int_{B_{r_0}}G_j(|\nabla
u_j|)\text{d}x\rightarrow \int_{B_{r_0}}G_\infty(|\nabla
u_0|)\text{d}x,
\end{align*}
and
  \begin{align*}\int_{B_{r_0}}G_j(|\nabla
u_j+\nabla \varphi|)\text{d}x\rightarrow \int_{B_{r_0}}G_\infty(|\nabla
u_0+\nabla \varphi|)\text{d}x.
\end{align*}
Then we obtain by \eqref{+21}, \eqref{+22'}, \eqref{+22}, \eqref{+23}, and $\lambda_{+j}\rightarrow0$
 \begin{align*}
 \int_{B_{r_0}}G_\infty(|\nabla
u_0)\text{d}x\leq \int_{B_{r_0}}G_\infty(|\nabla
u_0+\nabla \varphi|)\text{d}x.
\end{align*}
This implies that $u_0$ is a $G_\infty-$harmonic function in $ B_{r_0}$.

 Since $u_j\geq 0$ in $B_{1}$, $u_0\geq 0$ in $B_{r_0}$. Note that $u_0(0)=0$. The Harnack's inequality implies $u_0\equiv 0$ in $B_{r_0}$. Finally we get $u_0\equiv 0$ in $B_{1}$ due to the continuity of $u$ and the arbitrariness of $r_0$, which is a contradiction to \eqref{14413}.\hfill $\blacksquare$
 \end{pf*}
{}{ \begin{remark}
 We point out that the result in Lemma \ref{Harnack's inequality} is a special case of \cite[ Corollary 1.4]{L}, which has a slight different version of Harnack's inequality if $h\in L^m(\Omega)$ with $m>N$ (see proofs of \cite[Section 3 and 4]{L}). Therefore, one may obtain the Local Lipschitz continuity of nonnegative minimizers of $ \mathcal {J}(u)$ with $\lambda_+>0$, provided $h\in L^m(\Omega)$ with $m>N$.
 \end{remark}}


\begin{thebibliography}{99}

\bibitem{AF} Adams,~R.~A., Fournier,~J.~J.~F.: \emph{Sobolev Spaces}. Pure and Applied Mathematics. Amsterdam: Acad. Press (140) 2003.

\bibitem{ACF} Alt, H. M., Caffarelli, L. A., Friedman, A.: A free boundary problem for quasilinear elliptic
equations. \emph{Ann. Sc. Norm. Super. Pisa Cl. Sci.}, (4)11(1)(1984), 1-44.

\bibitem{AP} Alt, H. M., Phillips, D.: A free boundary problem for semilinear elliptic equations. \emph{J. Reine
Angew. Math.}, 368(1986), 63-107.

\bibitem{B} Braga, J. E. M.: On the Lipschitz regularity and asymptotic behaviour of
the free boundary for classes of minima of inhomogeneous
two-phase Alt-Caffarelli functionals in Orlicz spaces. \emph{Ann. Mat. Pura Appl.}, https://doi.org/10.1007/s10231-018-0755-7.


\bibitem{BM} Braga, J. E. M., Moreira, D. R.:  Uniform Lipschitz regularity for classes of minimizers in two phase
free boundary problems in Orlicz spaces with small density
on the negative phase. \emph{Ann. I. H. Poincar\'{e}- AN.}, 31(4)(2014), 823--850.

\bibitem{CL} Challal, S., Lyaghfouri, A.: Porosity of free boundaries in $A$-obstacle problems. \emph{Non. Anal.: TMA.}, 70(7)(2009), 2772-2778.

\bibitem{CLR} Challal, S., Lyaghfouri, A., Rodrigues, J. F.: On the $A$-obstacle problem and the Hausdorff measure of its free boundary. \emph{Ann. Mat. Pura Appl.}, 191(2012), 113-165.

    \bibitem{G}Gossez, J.: Nonlinear elliptic boundary value problems for equations with rapidly (or slowly) increasing coefficients. \emph{Trans. Am. Math. Soc.}, 190(1974), 163-205.


   \bibitem{DP} Danielli, D., Petrosyan, A.: A minimum problem with free boundary for a degenerate quasilinear
operator. \emph{Calc. Var. Partial Differ. Equ.}, 23(1)(2005), 97-124.


\bibitem{GG} Giaquinta, M., Giusti, E.: Differentiability of minima of non-differentiable functionals. \emph{Invent.
Math.}, 72(1983), 285-298.

\bibitem{KKPS} Karp, L., Kilpel\"{a}inen, T., Petrosyan, A., Shahgholian, H.: On the porosity of free boundaries
in degenerate variational inequalities. \emph{J. Differ. Equ.}, 164(2000), 110-117.



\bibitem{LU} Ladyzhenskaya, O. A., Ural'tseva, N. N.: \emph{Linear and Quasilinear Elliptic Equations}. vol. 46. Academic Press, New York, 1968.

 \bibitem{LS} Lee, K., Shahgholian, H.: Hausdorff measure and stability for the $p$-obstacle problem $(2 <
p <1)$. \emph{J. Differ. Equ.}, 195(2003), 14 -- 24.

\bibitem{LdT} Leit$\tilde{a}$o, R., de Queiroz, O. S., Teixeira,
E. V.: Regularity for degenerate two-phase free boundary problems. \emph{Ann. Inst. H. Poincar\'{e} Anal. Non Lineare.}, 32(4)(2015), 741 -- 762.

\bibitem{LT} Leit$\tilde{a}$o, R., Teixeira,
E. V.: Regularity and geometric estimates for minima of
discontinuous functionals. \emph{Rev. Mat. Iberoam.}, 31(1)(2015), 69 -- 108.

 \bibitem{L} Lieberman, G. M.: The natural generalization of the natural conditions of Ladyzhensaya and Uraltseva for elliptic
equations. \emph{Comm. Partial Differ. Equ.}, 16
(2\text{\&}3) (1991), 311 -- 361.


\bibitem{MW} Mart\'{\i}nez, S., Wolanski, N.: A minimum problem with free boundary
in Orlicz spaces. \emph{Adv. Math.}, 218(6) (2008), 1914 -- 1971.

\bibitem{FS}Fusco, N.,  Sbordone,C.: Some remarks on the regularity of minima of anisotropic integrals, Comm. Partial Differential Equations 18 (1993) 153--167.


\bibitem{P1} Phillips, D.: A minimization problem and the regularity of solutions in the presence of a free
boundary. \emph{Indiana Univ. Math. J.}, 32(1983), 1 -- 17.

\bibitem{P2} Phillips, D.: Hausdorff measure estimates of a free boundary for a minimum problem. \emph{Comm.
Partial Differ. Equ.}, 8(1983), 1409 -- 1454.

\bibitem{ZFZ} Zheng, J., Feng, B., Zhao, P.: Regularity of minimizers in the two-phase free boundary problems in Orlicz-Sobolev
spaces. \emph{Z. Anal. Anwend.}, 36(1)(2017), 37 -- 47.

\bibitem{ZG} Zheng, J., Guo, X.:~Lyapunov-type inequalities for $\psi$ - Laplacian equations. chinaXiv:201805. 00171.

\bibitem{ZZZ} Zheng, J., Zhang, Z., Zhao, P.: A minimum problem with two-phase free
boundary in Orlicz spaces. \emph{Monatsh. Math.}, 172(3-4) (2013), 441 -- 475.



\end{thebibliography}
\end{document}